\documentclass[10pt]{article}

\usepackage[bottom=0.5in]{geometry}

\usepackage{amsmath}
\usepackage{amssymb}
\usepackage{graphicx}
\usepackage{amsthm}
\usepackage{caption}
\usepackage{ wasysym }
\usepackage{tikz-cd} 
\usepackage{ stmaryrd }
\usepackage{xcolor}
\usepackage{scrextend}
\usepackage{mathtools}
\usepackage{xcolor}
\usepackage{bbm}
\usepackage{tabu}
\usepackage{tcolorbox}
\usepackage{dsfont}
\usepackage{pgfplots}
\usepackage{mathtools}
\usepackage{ wasysym }
\usepgfplotslibrary{polar}
\usepackage{ragged2e}
\usepackage{changepage}

\makeatletter
\newtheorem*{rep@theorem}{\rep@title}
\newcommand{\newreptheorem}[2]{%
\newenvironment{rep#1}[1]{%
 \def\rep@title{#2 \ref{##1}}%
 \begin{rep@theorem}}%
 {\end{rep@theorem}}}
\makeatother
\newreptheorem{theorem}{Theorem}
\newreptheorem{corollary}{Corollary}

\usepackage{cite}

\usepackage{float}
\usepackage{hyperref}

\usepackage[font=small,labelfont=bf]{caption}
\DeclareCaptionFont{small}{\small}

\makeatletter
\DeclareRobustCommand\widecheck[1]{{\mathpalette\@widecheck{#1}}}
\def\@widecheck#1#2{%
    \setbox\z@\hbox{\m@th$#1#2$}%
    \setbox\tw@\hbox{\m@th$#1%
       \widehat{%
          \vrule\@width\z@\@height\ht\z@
          \vrule\@height\z@\@width\wd\z@}$}%
    \dp\tw@-\ht\z@
    \@tempdima\ht\z@ \advance\@tempdima2\ht\tw@ \divide\@tempdima\thr@@
    \setbox\tw@\hbox{%
       \raise\@tempdima\hbox{\scalebox{1}[-1]{\lower\@tempdima\box
\tw@}}}%
    {\ooalign{\box\tw@ \cr \box\z@}}}
\makeatother

\def\centerarc[#1](#2)(#3:#4:#5);%
    {
    \draw[#1]([shift=(#3:#5)]#2) arc (#3:#4:#5);
    }
    
\NewEnviron{largeequation*}{%
    \begin{equation*}
    \scalebox{1.1}{$\BODY$}
    \end{equation*}
    }    
    
    \NewEnviron{smallequation*}{%
    \begin{equation*}
    \scalebox{0.8}{$\BODY$}
    \end{smallequation*}
    }    
    
\DeclareMathAlphabet\euscr{U}{eus}{m}{n}

\newcommand{\R}{\mathbb{R}}
\newcommand{\Z}{\mathbb{Z}}
\newcommand{\N}{\mathbb{N}}
\newcommand{\intersection}{\cap}

\newcommand{\isom}{\cong}
\newcommand{\surj}{\twoheadrightarrow}

\newcommand{\union}{\cup}
\newcommand{\Union}{\bigcup}
\newcommand{\Intersection}{\bigcap}

\newcommand{\inj}{\hookrightarrow}
\newcommand{\quotient}[2]{{\raisebox{.2em}{$#1$}\left/\raisebox{-.2em}{$#2$}\right.}}

\DeclarePairedDelimiterX{\inner}[2]{\langle}{\rangle}{#1, #2}

\newcommand{\C}{\mathbb{C}}
\newcommand{\ds}{\oplus}
\newcommand{\Ds}{\bigoplus}

\renewcommand{\o}{\circ}

\newcommand{\Q}{\mathbb{Q}}

\newcommand{\x}{\times}

\renewcommand{\x}{\otimes}

\renewcommand{\k}{\mathds{k}}

\newcommand{\twedge}{\vee}

\newcommand{\scrF}{\euscr{F}}

\newcommand{\Disj}{\bigsqcup}
\newcommand{\scrA}{\euscr{A}}

\renewcommand{\choose}[2]{\binom{#1}{#2}}

\newtheorem{theorem}{Theorem}[section]
\newtheorem{lemma}[theorem]{Lemma}
\newtheorem{proposition}[theorem]{Proposition}
\newtheorem{corollary}[theorem]{Corollary}

\newtheorem*{theorem*}{Theorem}

\theoremstyle{definition}
\newtheorem{definition}[theorem]{Definition}
\newtheorem{example}[theorem]{Example}
\newtheorem{observation}[theorem]{Observation}

\theoremstyle{remark}
\newtheorem{remark}[theorem]{Remark}

\DeclareMathOperator{\sep}{ | }

\DeclareMathOperator{\Hom}{Hom}

\DeclareMathOperator{\End}{End}
\DeclareMathOperator{\Span}{Span}

\DeclareMathOperator{\rk}{rk}

\DeclareMathOperator{\TL}{TL}

\DeclareMathOperator{\FI}{FI}
\DeclareMathOperator{\id}{id}
\DeclareMathOperator{\LS}{LS}

\DeclareMathOperator{\Comp}{Comp}
\DeclareMathOperator{\row}{row}
\DeclareMathOperator{\spn}{span}

\DeclarePairedDelimiter\floor{\lfloor}{\rfloor}

\definecolor{myblue}{rgb}{0.70,0.81,1}
\definecolor{myorange}{rgb}{1,0.72,0.4}
\definecolor{myred}{rgb}{1,0.58,0.58}
\definecolor{myyellow}{rgb}{1,0.88,0.56}
\definecolor{mylime}{rgb}{0.9,1,0}
\definecolor{mycyan}{rgb}{0.78,0.94,1}
\definecolor{mypink}{rgb}{0.78,0.94,1}
\definecolor{mymagenta}{rgb}{0.9,0.70,0.87}
\definecolor{myolive}{rgb}{0.83,0.87,0.54}
\definecolor{mypink}{rgb}{0.99,0.74,0.94}
\definecolor{mybrown}{rgb}{0.89,0.68,0.68}
 
 \newtcolorbox{custombox}[3][]
{
  colframe = #2!25,
  colback  = #2!10,
  coltitle = #2!20!black,  
  title    = #3,
  #1,
}
 
\usepackage{tikz}

\tikzset{
  pt/.style={insert path={node[scale=2]{.}}},
  dnup/.style={insert path={ [pt] .. controls +(0,1) and +(0,-1) .. +(#1,2) [pt]}},
  dndn/.style={insert path={ [pt] .. controls +(0,1) and +(0,1) .. +(#1,0) [pt]}},
  upup/.style={insert path={ [pt] .. controls +(0,-1) and +(0,-1) .. +(#1,0) [pt]}},
  sdndn/.style={insert path={ [pt] .. controls +(0,0.5) and +(0,0.5) .. +(#1,0) [pt]}},
  supup/.style={insert path={ [pt] .. controls +(0,-0.5) and +(0,-0.5) .. +(#1,0) [pt]}},
  bdndn/.style={insert path={ [pt] .. controls +(0,2) and +(0,2) .. +(#1,0) [pt]}},
  bupup/.style={insert path={ [pt] .. controls +(0,-2) and +(0,-2) .. +(#1,0) [pt]}},
  bbupup/.style={insert path={ [pt] .. controls +(0,-3) and +(0,-3) .. +(#1,0) [pt]}},
}

\renewcommand{\restriction}{\mathord{\upharpoonright}}

\providecommand{\keywords}[1]
{
  \small	
  \emph{Keywords:} #1
}

\providecommand{\amsclassification}[1]
{
  \small	
  \emph{MSC:} #1
}

\providecommand{\cis}[1]
{
  \small	
  \emph{Competing Interest Statement:} #1
}

\begin{document}
 
\title{ Topological actions of Temperley-Lieb monoids and representation stability}
\date{}
\author{Maithreya Sitaraman$^*$}
\date{\small \flushleft $^*$ Department of Mathematics, Columbia University, 2990 Broadway, New York, NY, United States \\ }
\maketitle

\begin{abstract}
We consider the Temperley-Lieb algebras $\TL_n(\delta)$ at $\delta = 1$. Since $\delta = 1$, we can consider the multiplicative monoid structure and ask how this monoid acts on topological spaces. Given a monoid action on a topological space, we get an algebra action on each homology group. The main theorem of this paper explicitly deduces the representation structure of the homology groups in terms of a natural filtration associated with our $\TL_n$-space. As a corollary of this result, we are able to study stability phenomena. There is a natural way to define representation stability in the context of $\TL_n(1)$, and the presence of filtrations enables us to define a notion of topological stability. We are able to deduce that a filtration-stable sequence of $\TL_n$-spaces results in representation-stable sequence of homology groups. This can be thought of as the analogue of the statement that the homology of configuration spaces forms a finitely generated $\FI$-module.
\end{abstract}

\amsclassification{81R05 $\cdot$ 16P70}


\cis{Declarations of interest: none.}

\let\thefootnote\relax\footnotetext{Email address: maithreya@math.columbia.edu}

{\small
\tableofcontents
}

\

\section{Introduction} 

\

Topological actions and representation stability go hand in hand, and it therefore makes sense to study both together. Indeed, the homology of natural families of spaces provides rich examples of representation stability. The foundational paper on representation stability \cite{rs1} (Church, Ellenberg, Farb) shows, for instance, that the homology of configuration spaces is representation-stable. Temperley-Lieb algebras naturally form a chain of inclusions. Even the earlier work on Temperley-Lieb algebras before diagrammatics became fashionable recognized the importance of inductive procedures and constructions, making use of the fact that they formed a chain \cite{tl5} (Goodman, Wenzl). Many modern treatments \cite{tl1}, \cite{tl2}, \cite{tl6} (Ridout, Saint-Aubin, Belletete) tend to involve diagrammatics, and the inclusions now have the appealing diagrammatic interpretation of adding loose strands. Moreover, this process of adding loose strands carries over to the representation theory, at least to the standard modules (see \cite{tl1} (Ridout, Saint-Aubin) for terminology). Needless to say, induction and restriction on this chain play a crucial role in understanding the representation theory, much like the analogue of symmetric groups. All this is suggestive that there ought to be a natural notion of representation stability for Temperley-Lieb algebras. Since topological actions and representation stability go hand in hand, it is therefore natural to study topological actions of Temperley-Lieb algebras.

\

We are interested in topological actions of the Temperley-Lieb monoids, that is: 

\begin{definition}\label{deftln}
The Temperley-Lieb monoid, $(\TL_n, \cdot)$ is the monoid of Temperley-Lieb diagrams from $n$ to $n$ points, with multiplication being concatenation, and with circles evaluated to $1$. 
\end{definition}

\begin{remark}
Note that for any field $\k$, $\k[\TL_n]$ is then the specialization of the Temperley-Lieb algebra over $\k$, $\TL_n(\delta)$, at $\delta = 1$.
\end{remark}

\begin{definition}
Given a topological space $X$, a topological action of $\TL_n$ is a map $\TL_n \to \Hom(X,X) = C(X,X)$, i.e a map to the collection of continuous functions on $X$.
\end{definition}

\

To our knowledge, topological actions of $\TL_n$ have not been studied previously. The representation stability of Temperley-Lieb algebras has debatably been considered before, albeit without the author mentioning representation stability. Indeed, in \cite{tl3} (Moore), the author considers representations of $\TL(\infty)$, and classifies ``link state representations'' which are indecomposable or irreducible. However, $\TL(\infty)$ is likely to be of ``wild type'', just as the infinite symmetric group $S(\infty)$ is of wild type. For $S(\infty)$, the attention therefore was shifted to studying sufficiently nice families of representations, such as tame representations, factor representations, and their generalizations (see e.g \cite{si1} (Okounkov) ). Extracting a notion of a stable representation from an infinite algebra $A(\infty)$ sometimes does not produce the notion of stability we want in representation stability, see for example \cite{si4} (Vershik, Nessonov). The author of \cite{tl3} (Moore) does not make any suggestion regarding what families of $\TL(\infty)$-representations could be considered stable. We therefore take a different approach to define a notion of stability, and our view is that it is the natural analogue to the definition of $\FI$-modules and hence much of the intuition from the theory of $\FI$-modules can be carried over. For our purposes, this analogy with $\FI$-modules is important, since the statement that filtration-stability of $\TL_n$-spaces implies representation stability (Theorem \ref{fsirs}) can be thought of as the analogue of the statement that the homology of configuration spaces forms a finitely generated $\FI$-module, as in \cite{rs1} (Church, Ellenberg, Farb). There have been other stability phenomena which have been studied in the context of Temperley-Lieb algebras, such as their homology and cohomology groups \cite{tl4} (Boyd, Hepworth) - of course, these have a very different flavor from that of representation stability. 

\

To our knowledge, Temperley-Lieb algebras have not been studied in the representation stability literature, or within the broader context of representation stability and $\FI$-modules. It appears that much of the work in representation stability has focussed on algebraic objects which are either close to symmetric groups \cite{rs2} \cite{rs3} \cite{rs6} (Wilson, Putman, Sam, Gunturkun, Snowden ) or are close to Lie groups \cite{rs3} \cite{rs5} (Sam, Snowden, Putman). Diagrammatically defined chains of algebras appear to have not been considered as objects whose representation category can be studied through the lens of representation stability. Diagrammatics and representation stability have, however, been uttered in the same breadth, but in a different sense: in \cite{dc1} (Barter, Entova-Aizenbud, Heidersdorf) the authors produce a functor from the category of $\FI$-modules modulo finite length $\FI$-modules to the abelian envelope of the Deligne category. The chain with respect to which one is considering representation stability there is of course still the chain of symmetric groups. Thus, representation stability with respect to a chain of diagrammatically defined algebras is not considered in \cite{dc1} (Barter, Entova-Aizenbud, Heidersdorf). 

\

Deducing the representation structure of homology groups given a topological action often does not have a known explicit formula. For instance decomposing homology groups of configuration spaces into irreducible representations of $S_n$ remains open even for simple manifolds \cite{ta1} (Hersh, Reiner). Sometimes, however, an explicit formula exists, is classically understood, and is relatively simple (though usually non-trivial to prove), as is the case for various Chevalley-Weil formulae, for example \cite{ta2} \cite{ta3} (Kloosterman, Grunewald, Larsen, Lubotzky, Malestein). Our main theorem shows us that topological actions of Temperley-Lieb algebras lie somewhere in the middle of this spectrum: on one hand, the formula we obtain for the representation structure of homology is explicit, indicating a degree of rigidity that $\TL_n$-spaces exhibit. On the other hand, our formula that describes how homology groups decompose is quite interesting, and has a nice combinatorial formulation. 

\

Given a suitable topological space $X$ with a suitable topological action $\TL_n \to \Hom(X,X)$, the main result of the paper (Theorem \ref{main}) deduces the structure of homology $H_*(X)$ as a representation of $\TL_n$. It turns out that under the hypotheses of the theorem, the representation structure of homology can be read off from an auxiliary construction $\scrF$ (see Definition \ref{fr}). Noting that $\TL_n$ is not semisimple, the reader may find our result (Theorem \ref{main}) somewhat surprising. However, to reconcile non-semisimplicity with our result, the reader may note that the result fails if assumptions on $X$ are relaxed.

\

As a corollary of the main result of the paper (Corollary \ref{fsirs}), we deduce a result which relates topological stability to representation stability in the context of $\TL_n$. If we think of configuration spaces as the notion of topological stability for symmetric groups, then a celebrated theorem in the foundational paper on $\FI$ modules (i.e \cite{rs1} (Church, Ellenberg, Farb), Section 6) shows that the homology of configuration spaces is a finitely generated $\FI$-module. Our theorem can be thought of as the analogue of this statement for the chain of Temperley-Lieb monoids.

\

\section{Quick review: Temperley-Lieb algebras and their representation theory}

\

\subsection{The definition of Temperley-Lieb algebras}

Given a field $\k$ and $\delta \in \k$, the Temperley-Lieb algebra $\TL_n(\delta)$ is the algebra of all crossingless matchings from $n$ points to $n$ points, where multiplication is concatenation, circles are evaluated to $\delta$, and addition is formal. In the current paper, we focus on the case when $\delta = 1$, and we are not interested in the formal additive structure of $\TL_n(1)$ but only on the multiplicative monoid of  planar diagrams that generates $\TL_n(1)$, which we refer to as simply $\TL_n$, as in Definition \ref{deftln}.

\

For example, for a field $\k$ and $\delta \in \k$, $\TL_3(\delta)$ is the algebra comprising of formal $\k$-linear combinations of the elements: $\begin{tikzpicture}[scale = 0.2, baseline={(0,0.1)}]
    \draw (0,0) -- (0,2);
    \draw[fill = black] (0,0) circle(0.15);
    \draw[fill = black] (0,2) circle(0.15);
    \draw (1,0) -- (1,2);
    \draw[fill = black] (1,0) circle(0.15);
    \draw[fill = black] (1,2) circle(0.15);   
    \draw (2,0) -- (2,2);
    \draw[fill = black] (2,0) circle(0.15);
    \draw[fill = black] (2,2) circle(0.15);     
  \end{tikzpicture}$ , $\begin{tikzpicture}[scale = 0.2, baseline={(0,0.1)}]
    \draw (1,2) [upup=1];
    \draw (1,0) [dndn=1];
    \draw (0,0) -- (0,2);
    \draw[fill = black] (0,0) circle(0.15);
    \draw[fill = black] (0,2) circle(0.15);
  \end{tikzpicture}$, $\begin{tikzpicture}[scale = 0.2, baseline={(0,0.1)}]
    \draw (0,2) [upup=1];
    \draw (0,0) [dndn=1];
    \draw (2,0) -- (2,2);
    \draw[fill = black] (2,0) circle(0.15);
    \draw[fill = black] (2,2) circle(0.15);
  \end{tikzpicture}$ , $\begin{tikzpicture}[scale = 0.2, baseline={(0,0.1)}]
    \draw (1,2) [upup=1];
    \draw (0,0) [dndn=1];
    \draw (2,0) [dnup=-2];
  \end{tikzpicture}$ , $\begin{tikzpicture}[scale = 0.2, baseline={(0,0.1)}]
    \draw (0,2) [upup=1];
    \draw (1,0) [dndn=1];
    \draw (0,0) [dnup=2];
  \end{tikzpicture}$. Multiplication is given by, for example: 
  
\

$$\begin{tikzpicture}[scale = 0.2, baseline={(0,0.1)}]
    \draw (1,2) [upup=1];
    \draw (0,0) [dndn=1];
    \draw (2,0) [dnup=-2];
  \end{tikzpicture} \cdot \begin{tikzpicture}[scale = 0.2, baseline={(0,0.1)}]
    \draw (1,2) [upup=1];
    \draw (1,0) [dndn=1];
    \draw (0,0) -- (0,2);
    \draw[fill = black] (0,0) circle(0.15);
    \draw[fill = black] (0,2) circle(0.15);
  \end{tikzpicture}  = \begin{tikzpicture}[scale = 0.2, baseline={(0,0.1)}]
    \draw (1,2) [upup=1];
    \draw (0,0) [dndn=1];
    \draw (2,0) [dnup=-2];
     \draw (1,0) [upup=1];
    \draw (1,-2) [dndn=1];
    \draw (0,-2) -- (0,0);
    \draw[fill = black] (0,-2) circle(0.15);
    \draw[fill = black] (0,0) circle(0.15);   
  \end{tikzpicture} = \begin{tikzpicture}[scale = 0.2, baseline={(0,0.1)}]
    \draw (1,2) [upup=1];
    \draw (1,0) [dndn=1];
    \draw (0,0) -- (0,2);
    \draw[fill = black] (0,0) circle(0.15);
    \draw[fill = black] (0,2) circle(0.15);
  \end{tikzpicture} $$

and 

$$\begin{tikzpicture}[scale = 0.2, baseline={(0,0.1)}]
    \draw (1,2) [upup=1];
    \draw (1,0) [dndn=1];
    \draw (0,0) -- (0,2);
    \draw[fill = black] (0,0) circle(0.15);
    \draw[fill = black] (0,2) circle(0.15);
  \end{tikzpicture} \cdot \begin{tikzpicture}[scale = 0.2, baseline={(0,0.1)}]
    \draw (1,2) [upup=1];
    \draw (1,0) [dndn=1];
    \draw (0,0) -- (0,2);
    \draw[fill = black] (0,0) circle(0.15);
    \draw[fill = black] (0,2) circle(0.15);
  \end{tikzpicture}  = \begin{tikzpicture}[scale = 0.2, baseline={(0,0.1)}]
    \draw (1,2) [upup=1];
    \draw (1,0) [dndn=1];
    \draw (0,0) -- (0,2);
    \draw[fill = black] (0,0) circle(0.15);
    \draw[fill = black] (0,2) circle(0.15);
     \draw (1,0) [upup=1];
    \draw (1,-2) [dndn=1];
    \draw (0,-2) -- (0,0);
    \draw[fill = black] (0,-2) circle(0.15);
    \draw[fill = black] (0,0) circle(0.15);   
  \end{tikzpicture} = \delta \begin{tikzpicture}[scale = 0.2, baseline={(0,0.1)}]
    \draw (1,2) [upup=1];
    \draw (1,0) [dndn=1];
    \draw (0,0) -- (0,2);
    \draw[fill = black] (0,0) circle(0.15);
    \draw[fill = black] (0,2) circle(0.15);
  \end{tikzpicture} $$
  
 \
 
 $\TL_n(\delta)$ admits the very useful presentation (see e.g \cite{tl1}) with $n-1$ generators $\{u_1,u_2,...,u_{n-1}\}$ subject to the relations: 
 
 \
 
\begin{equation*}
\begin{cases}
u_i^2 = \delta \cdot u_i \textrm{ for all $i \in \{1,2,...,n-1\}$ } & \textrm{ (the idempotent relation)}\\ 
u_i u_{i+1} u_i = u_i \textrm{ for all $i \in \{1,2,...,n-2\}$ } & \textrm{ (the upper neighbor relation)}\\ 
u_i u_{i-1} u_i = u_i \textrm{ for all $i \in \{2,...,n-1\}$ } & \textrm{ (the lower neighbor relation)}\\ 
u_i u_j = u_j u_i \textrm{ whenever $|i-j| \ge 2$ } & \textrm{ (the long-distance relation)} \\ 
\end{cases}
\end{equation*}

\

Here $u_i$ can be represented diagrammatically as the diagram which comprises of a single cup and cap at position $i$, and straight lines elsewhere. For example, $u_2 \in \TL_5$ is given by the diagram $\begin{tikzpicture}[scale = 0.2, baseline={(0,0.1)}]
    \draw (1,2) [upup=1];
    \draw (1,0) [dndn=1];
    \draw (3,0) -- (3,2);
    \draw[fill = black] (3,0) circle(0.15);
    \draw[fill = black] (3,2) circle(0.15);
    \draw (4,0) -- (4,2);
    \draw[fill = black] (4,0) circle(0.15);
    \draw[fill = black] (4,2) circle(0.15);
    \draw (0,0) -- (0,2);
    \draw[fill = black] (0,0) circle(0.15);
    \draw[fill = black] (0,2) circle(0.15);
  \end{tikzpicture}$.
 
\

\subsection{Link state representations of Temperley-Lieb algebras}

\

Link state representations are an important class of representations of $\TL_n(\delta)$. When $\TL_n(\delta)$ is semisimple (this is true when, for example, $\k = \mathbb{C}$ and $\delta \not= q+q^{-1}$ for $q$ a root of unity, see for example \cite{tl1}), then every representation of $\TL_n(\delta)$ is a direct sum of link-state representations. When $\delta$ is not semisimple, there are whole host of other representations, but we will not talk about these in this Subsection. The reader may refer to \cite{tl2} (Belletete, Ridout, Saint-Aubin) for a classification of indecomposables when $\k = \C$.

\

For $p$ an integer $\le \frac{n}{2}$, an $(n,p)$ link state is a collection $p$ cups and $n-2p$ lines. For example, the $(6,2)$ link states are: $\left(\begin{tikzpicture}[scale = 0.2, baseline={(0,0.1)}]
    \draw (0,0) [upup=1];
    \draw (2,0) [upup=1];
    \draw (4,0) -- (4,-1);
    \draw[fill = black] (4,0) circle(0.15); 
    \draw (5,0) -- (5,-1);
    \draw[fill = black] (5,0) circle(0.15); 
 \end{tikzpicture}\right)$, $\left(\begin{tikzpicture}[scale = 0.2, baseline={(0,0.1)}]
    \draw (0,0) [upup=1];
    \draw (2,0) -- (2,-1);
    \draw[fill = black] (2,0) circle(0.15);     
    \draw (3,0) [upup=1];
    \draw (5,0) -- (5,-1);
    \draw[fill = black] (5,0) circle(0.15); 
 \end{tikzpicture}\right)$, $\left(\begin{tikzpicture}[scale = 0.2, baseline={(0,0.1)}]
     \draw (0,0) -- (0,-1);
    \draw[fill = black] (0,0) circle(0.15);      
    \draw (1,0) [upup=1];
    \draw (3,0) [upup=1];
    \draw (5,0) -- (5,-1);
    \draw[fill = black] (5,0) circle(0.15); 
 \end{tikzpicture}\right)$, $\left(\begin{tikzpicture}[scale = 0.2, baseline={(0,0.1)}]
    \draw (0,0) [upup=1];
    \draw (2,0) -- (2,-1);
    \draw[fill = black] (2,0) circle(0.15);
     \draw (3,0) -- (3,-1);
    \draw[fill = black] (3,0) circle(0.15);      
    \draw (4,0) [upup=1];
 \end{tikzpicture}\right)$, $\left(\begin{tikzpicture}[scale = 0.2, baseline={(0,0.1)}]
    \draw (0,0) -- (0,-1);
    \draw[fill = black] (0,0) circle(0.15);
    \draw (1,0) [upup=1];    
     \draw (3,0) -- (3,-1);
    \draw[fill = black] (3,0) circle(0.15);      
    \draw (4,0) [upup=1];
 \end{tikzpicture}\right)$, $\left(\begin{tikzpicture}[scale = 0.2, baseline={(0,0.1)}]
    \draw (0,0) -- (0,-1);
    \draw[fill = black] (0,0) circle(0.15);
     \draw (1,0) -- (1,-1);    
     \draw[fill = black] (1,0) circle(0.15);      
    \draw (2,0) [upup=1];         
    \draw (4,0) [upup=1];
 \end{tikzpicture}\right)$, $\left(\begin{tikzpicture}[scale = 0.2, baseline={(0,0.1)}]
    \draw (0,0) [upup=3];
    \draw (1,0) [supup=1];
    \draw (4,0) -- (4,-1);
    \draw[fill = black] (4,0) circle(0.15); 
    \draw (5,0) -- (5,-1);
    \draw[fill = black] (5,0) circle(0.15); 
 \end{tikzpicture}\right)$, $\left(\begin{tikzpicture}[scale = 0.2, baseline={(0,0.1)}]
     \draw (0,0) -- (0,-1);
    \draw[fill = black] (0,0) circle(0.15);      
    \draw (1,0) [upup=3];
    \draw (2,0) [supup=1];
    \draw (5,0) -- (5,-1);
    \draw[fill = black] (5,0) circle(0.15); 
 \end{tikzpicture}\right)$, $\left(\begin{tikzpicture}[scale = 0.2, baseline={(0,0.1)}]
    \draw (0,0) -- (0,-1);
    \draw[fill = black] (0,0) circle(0.15);
     \draw (1,0) -- (1,-1);    
     \draw[fill = black] (1,0) circle(0.15);      
    \draw (2,0) [upup=3];         
    \draw (3,0) [supup=1];
 \end{tikzpicture}\right)$.

\

\

Let $M_n$ be the set $\Union_{p=1}^{\floor{\frac{n}{2}}} \{(n,p) \textrm{ link states}\}$. Abusing notation slightly, let us also let $M_n$ denote the $\k$-linear span of that set. $M_n$ admits an action of $\TL_n(\delta)$ by concatenating a diagram of $\TL_n(\delta)$ with a link state and evaluated any circles to $\delta$. If $\delta = 1$, this becomes a well defined action on the set $M_n$. For the moment, we take $\delta$ arbitrary and present an example to illustrate the action: 

\

$$\begin{tikzpicture}[scale = 0.2, baseline={(0,0.1)}]
    \draw (1,2) [upup=1];
    \draw (1,0) [dndn=1];
    \draw (0,0) -- (0,2);
    \draw[fill = black] (0,0) circle(0.15);
    \draw[fill = black] (0,2) circle(0.15);
     \draw (3,2) [upup=3];
    \draw (4,2) [supup=1];
     \draw (3,0) [dndn=1];
    \draw (5,0) [dndn=1];      
  \end{tikzpicture} \cdot \left(\begin{tikzpicture}[scale = 0.2, baseline={(0,0.1)}]
    \draw (0,0) -- (0,-1);
    \draw[fill = black] (0,0) circle(0.15);
    \draw (1,0) [upup=1];    
     \draw (3,0) -- (3,-1);
    \draw[fill = black] (3,0) circle(0.15);      
    \draw (4,0) [upup=1];
 \end{tikzpicture}\right) = \begin{tikzpicture}[scale = 0.2, baseline={(0,0.1)}]
    \draw (1,2) [upup=1];
    \draw (1,0) [dndn=1];
    \draw (0,0) -- (0,2);
    \draw[fill = black] (0,0) circle(0.15);
    \draw[fill = black] (0,2) circle(0.15);
     \draw (3,2) [upup=3];
    \draw (4,2) [supup=1];
     \draw (3,0) [dndn=1];
    \draw (5,0) [dndn=1];      
    \draw (0,0) -- (0,-1);
    \draw[fill = black] (0,0) circle(0.15);
    \draw (1,0) [upup=1];    
     \draw (3,0) -- (3,-1);
    \draw[fill = black] (3,0) circle(0.15);      
    \draw (4,0) [upup=1];    
  \end{tikzpicture} = \delta \left(\begin{tikzpicture}[scale = 0.2, baseline={(0,0.1)}]
    \draw (0,0) -- (0,-1);
    \draw[fill = black] (0,0) circle(0.15);
    \draw (1,0) [upup=1];    
    \draw (3,0) [upup=3];         
    \draw (4,0) [supup=1];
 \end{tikzpicture}\right) $$

\

$M_n$ admits a natural filtration as follows: the action of a Temperley-Lieb element can never decrease the number of cups of a link state. For example, in the above example, the action increases the number of cups from $2$ to $3$. So, let $M_{n,p} = \Union_{q \ge p}^{\floor{\frac{n}{2}}} \{(n,q) \textrm{ link states}\}$. 

\

Then, the $\TL_n$ action respects the filtration

$$M_{n,0} \supset M_{n,1} \supset M_{n,2} \supset ... \supset M_{n, \floor{\choose{n}{2}}}$$

\

Therefore, it is natural to define the representations $V_{n,p} = \quotient{M_{n,p}}{M_{n,p+1}}$ when $p < \floor{\choose{n}{2}}$ and $V_{n,\floor{\choose{n}{2}}} = M_{n, \floor{\choose{n}{2}}}$. $\{V_{n,p}\}$ are then well defined representations of $\TL_n(\delta)$. The representation $V_{n,p}$ is called the standard representation of $\TL_n(\delta)$ of $(n,p)$ link states. When $\k = \C$ and $\TL_n(\delta)$ is semisimple, $\{V_{n,p}\}_{p \in \{0,1,2,...,\floor{\choose{n}{2}}\}}$ form a complete set of irreducible representations of $\TL_n(\delta)$ (see for example \cite{tl1}). 

\

\section{Describing topological actions of Temperley-Lieb algebras} \label{sdtatla}

\

Let $(\TL_n, \cdot)$ be the Temperley-Lieb monoid, as in Definition \ref{deftln}. A topological action of $\TL_n$ on a topological space $X$ is a map

$$(\TL_n, \cdot) \mapsto (\Hom(X,X), \circ)$$

where $(\Hom(X,X), \circ)$ is of course the monoid of continuous maps from $X$ to itself, under composition. It is important to note that given a topological action $(\TL_n, \cdot) \mapsto (\Hom(X,X), \circ)$, the induced action on (co)homology is an algebra representation, i.e given a field $\k$ and considering $\TL_n(1)$ as an algebra over $\k$, we get a map $\TL_n(1) \to \End(H_k(X, \k))$ for each $k \in \N$. We will henceforth suppress the field $\k$ in our notation for homology groups.

\

A topological space $X$ that admits such a map $(\TL_n, \cdot) \mapsto (\Hom(X,X), \circ)$ will be called a \underline{$\TL_n$-space}.

\

\subsection{Topological translation of the Temperley-Lieb relations}

\

The first natural step to study such an action would be to look at the relations of $(\TL_n, \cdot)$ (with respect to a natural presentation), and translate what these relations mean in topological language. We will use the most common presentation of  $(\TL_n, \cdot)$, namely generated by $u_1,...,u_{n-1}$ subject to the relations:  

\begin{equation*}
\begin{cases}
u_i^2 = u_i \textrm{ for all $i \in \{1,2,...,n-1\}$ } & \textrm{ (the $\delta$-idempotent relation)}\\ 
u_i u_{i+1} u_i = u_i \textrm{ for all $i \in \{1,2,...,n-2\}$ } & \textrm{ (the upper neighbor relation)}\\ 
u_i u_{i-1} u_i = u_i \textrm{ for all $i \in \{2,...,n-1\}$ } & \textrm{ (the lower neighbor relation)}\\ 
u_i u_j = u_j u_i \textrm{ whenever $|i-j| \ge 2$ } & \textrm{ (the long-distance relation)} \\ 
\end{cases}
\end{equation*}

\

We now will translate each of the three relations into topological language: 

\

\begin{lemma}[Translating the idempotent relation] \label{tir} 
Each $u_i$ maps to a retraction map. 
\end{lemma} 

\begin{proof}
The idempotents in $(\Hom(X,X), \circ)$ are precisely retractions. 
\end{proof}

\

We will denote the retraction map associated to $u_i$ by $r_i$, and the subspace onto which $r_i$ retracts by $A_i$.

\

\begin{lemma}[Translating the neighbor relation] \label{tnr} 
$r_{i\pm1} \restriction_{A_i}: A_i \to A_{i\pm1}$ and $r_{i} \restriction_{A_{i\pm1}}: A_{i\pm1} \to A_{i}$ are homeomorphisms. Moreover, these two homeomorphisms are inverses of each other. 
\end{lemma} 

\begin{proof}
The neighbor relation tells us that $u_i u_{i\pm1} u_i = u_i$, and this therefore, $r_i r_{i\pm1} r_i = r_i$. Since $r_i(x) = x$ for every $x \in A_i$, we conclude that for every $x \in A_i$ $r_i r_{i\pm1}(x) = x$. That is, $r_i r_{i\pm1} = \id_{A_i}$. This tells us that $r_i \restriction_{A_{i\pm1}}$ and $r_{i\pm1} \restriction_{A_i}$ are inverses of each other, and thus, in particular they are both homeomorphisms.
\end{proof}

\

\begin{remark} \label{tnrrk}
It is natural to wonder whether Lemma \ref{tnr} is equivalent to the neighbor relation, or whether it is a weaker statement. Even though we only needed to take $x \in A_i$, it turns out that Lemma \ref{tnr} is in fact equivalent to the neighbor relation. The reason is that assuming the Lemma \ref{tnr} holds, and given any $x \in X$ arbitrary, $r_i(x) \in A_i$ and thus, applying the Lemma to $r_i(x)$, we deduce that $r_i r_{i\pm1} r_i (x) = r_i(x)$. Since $x$ was arbitrary, we deduce the Temperley-Lieb relation. This justifies the fact that we called it a translation.
\end{remark}

\

\begin{observation}[Composition of commuting retractions] \label{compcommrretr}
Let $X$ be a topological space, and let  $r: X \to A$ and $r': X \to A'$ be commuting retractions, i.e $r \circ r' = r' \circ r$. Then, $r \circ r'$ is a retraction onto $A \intersection A'$. 
\end{observation}

\begin{proof}
Given $x \in X$, $r r'(x) = r(r'(x)) \in A$ but also, by commutativity, $r r'(x) = r'(r(x)) \in A'$. Thus, $r r'(x) \in A \intersection A'$ for every $x \in X$. Moreover, if $x \in A \intersection A'$, then both $r$ and $r'$ fix $x$, and so $r r'(x) = r(x) = x$. Thus, $r r'$ is a retraction onto $A \intersection A'$. 
\end{proof}

\

\begin{lemma}[Translating the long-distance relation] \label{tdr} 
If $\{i_1,i_2,...,i_m\} \subset \{1,2,...,n-1\}$ is a set such that $|i_k - i_l| \ge 2$ for all $k \not= l$, then $r_{i_1} r_{i_2} ... r_{i_m}: X \to A_{i_1} \intersection ... \intersection A_{i_m}$ is a retraction onto $A_{i_1} \intersection ... \intersection A_{i_m}$.
\end{lemma} 

\begin{proof}
The long-distance relation tells us that $u_{i_k} u_{i_l} = u_{i_l} u_{i_k}$ for all $k \not= l$, and thus we have a family of commuting retractions. Applying Observation  \ref{compcommrretr} recursively, we see that $r_{i_1} r_{i_2} ... r_{i_m}$ is a retraction onto $A_{i_1} \intersection ... \intersection A_{i_m}$.
\end{proof}

\

At this point, the reader should be able to form some kind of sketchy picture in their heads for how topological actions of Temperley-Lieb algebras behave. Before we delve deeper, this would be a good point for the reader to take a look at some of the pictures in Subsection \ref{ssets} to improve this intuition. 

\

\subsection{Any neighbor intersection equals the full intersection} \label{sstfi}

\

 In this Subsection we will demonstrate an important fact, that any point in a neighbor intersection $A_i \intersection A_{i+1}$ must in fact lie in the full intersection  $A_1 \intersection A_2 \intersection ... \intersection A_{n-1}$. This will serve as one of the key foundational stones for the rest of the paper.

\

\begin{lemma}[neighbor intersection $=$ full intersection] \label{nifi} 
Suppose that $x \in A_i \intersection A_{i+1}$. Then, in fact, $x \in A_1 \intersection A_2 \intersection ... \intersection A_{n-1}$.
\end{lemma} 

\begin{proof}
By symmetry and induction, it suffices to show that $x \in A_{i+2}$. The trick now is to look at the ``cycle map''  $r_i r_{i+1} r_{i+2}$. 

\

\begin{figure}[H]
\centering
\begin{tikzpicture}[scale = 1.5]
\draw (-0.9,0) circle (1);
\draw (+0.9,0) circle (1);
\draw (0,1.55) circle (1);
\node at (-0.9,0) {$A_i$};
\node at (0.9,0) {$A_{i+1}$};
\node at (0,1.55) {$A_{i+2}$};
\draw[fill = black] (0,0) circle(0.05);
\draw[->, blue, thick] (0.9,-0.4) -- (-0.9,-0.4); 
\node at (0,-0.75) {\color{blue} $\mathbf{r_{i}}$};
\draw[->, blue, thick] (-1.185,0.285) -- (-0.285,1.835); 
\draw[->, blue, thick] (0.285,1.835) -- (1.185,0.285); 
\node at (0,-0.14) {$x$};
\node at (-1.2,1.1) {\color{blue} $\mathbf{r_{i+2}}$};
\node at (1.2,1.1) {\color{blue} $\mathbf{r_{i+1}}$};
\end{tikzpicture} 
\caption{A schematic of the cycle map}
\end{figure}

\

Since $x \in A_i$, it follows by Lemma \ref{tdr} that $r_{i+2}(x) \in A_i \intersection A_{i+2}$. Therefore, on one hand, by Lemma \ref{tnr} 

$$r_i r_{i+1} r_{i+2} x  = (r_i r_{i+1}) (r_{i+2}(x)) = \id_{A_i} r_{i+2} (x) = r_{i+2}(x) \in A_{i+2}$$

On the other hand, by associativity of multiplication, Lemma \ref{tnr}, and the fact that $x \in A_{i+1}$

$$r_i r_{i+1} r_{i+2} (x) = r_i (r_{i+1} r_{i+2})(x) = r_i \id_{A_{i+1}} x = r_i(x) = x$$

Therefore, comparing the two expressions, 

$$x= r_{i+2}(x) \in A_{i+2}$$
\end{proof}

\

Action of $\TL_n$ on a topological space $X$ induces an action on each homology group $H_k(X)$. For the reader's convenience, we will now state a homological version of Lemma \ref{nifi} which will be useful later in the paper. The reader may feel free to skip over this for now and return to this at a later point.

\begin{corollary}[homological version of neighborhood intersection = full intersection] \label{hnifi}
For each $A_j$, let $\iota: A_j \inj X$ denote the inclusion map. As before Let $r_j: X \to A_j \inj X$ denote the retraction map. Then, for $\alpha \in H_k(X)$, if there is some $j$ for which $\alpha \in \iota_* H_k(A_j)$ and $\alpha \in \iota_* H_k(A_{j+1})$, then in fact $(r_i)_* \alpha = \alpha$ for all $i$, i.e $\alpha$ is fixed by all $(r_i)_*$.
\end{corollary}

\begin{proof}
By Observation \ref{inclinjhom}, $\iota_*: H_k(A_j) \inj H_k(X)$ is an injection of homology groups. For each $r_i: X \to A_i \inj X$, let $(r_i)_*: H_k(X) \to H_k(A_i) \inj H_k(X)$ denote the corresponding map on homology. Since $\alpha \in \iota_* H_k(A_j), \iota_* H_k(A_{j+1})$, we have that $(r_j) _*\alpha = \alpha$ and $(r_{j+1})_* \alpha = \alpha$. By the same proof as Lemma \ref{nifi} (by the cycle map trick), we therefore have that $(r_i)_* \alpha = \alpha$ for all $i$.
\end{proof}

\

We conclude this subsection with an example of a $\TL_n$-space where $A_i \intersection A_{i+1} = \emptyset$ for all $i\in\{1,2,..,n-2\}$. It is due to examples like this that, later on, we will need to consider the image of the ideal $\overset{\sim}{\TL_n}$ in Theorem \ref{main}. 

\

\begin{example}[An example where $A_i \intersection A_{i+1} = \emptyset$ for all $i\in\{1,2,..,n-2\}$]
Let $X$ consist of $5$ lines, denoted by $A_1, A_2, A_3, A_4, A_5$, depicted below. With the action of $\TL_n$ described below, this constitutes an example of a $\TL_n$ space where $A_i \intersection A_{i+1} = \emptyset$ for all $i\in\{1,2,..,n-2\}$

\begin{figure}[H]
\centering
\begin{tikzpicture}
\draw (3,8) --(3,6);
\draw (3,6) --(1,4);
\draw (3,6) --(5,4);
\draw (1,4) --(3,2);
\draw (5,4) --(3,2);

   \node at (2.7,7) {$\mathbf{A_3}$};
   \node at (1.85,5.25) {$\mathbf{A_1}$};
   \node at (4.15,5.25) {$\mathbf{A_5}$};
   \node at (1.85,2.75) {$\mathbf{A_4}$};
   \node at (4.15,2.75) {$\mathbf{A_2}$};
   
   \draw[fill] (3,8) circle (3pt);
   \draw[fill] (3,6) circle (3pt);
   \draw[fill] (1,4) circle (3pt);
   \draw[fill] (5,4) circle (3pt);
   \draw[fill] (3,2) circle (3pt);

   \node at (3.2,8.2) {a};
   \node at (3.2,6.2) {b};
   \node at (0.8,4.2) {c};
   \node at (5.2,4.2) {d};
   \node at (3.2,1.8) {e};

\end{tikzpicture}
\caption{We define the action of $\TL_n$ on $X$ as follows: $r_1,r_2,r_3,r_4,r_5$ will be retractions onto $A_1,A_2,A_3,A_4,A_5$ in such a way that if $|i - j| \ge 2$, then $r_j(A_i)$ is a point and if $|i-j| = 1$, then $r_i(A_j) = A_i$. }
\end{figure}
We may explicitly define $\{r_i\}_{i=1}^5$ by considering actions on the points $a,b,c,d,e$.  
\begin{equation*}
\begin{split}
r_3(a) &= r_3(e) = a, r_3(b) = r_3(c) = r_3(d) = b \\
r_1(b) &= r_1(a) = r_1(d) = b, r_1(c) = r_1(e) = c \\
r_5(b) &= r_5(a) = r_5(c) = b, r_5(d) = r_5(e) = d \\
r_4(c) &= r_4(b) = c, r_4(e) = r_4(a) = r_4(d) = e \\
r_2(d) &= r_2(b) = d, r_2(e) = r_2(c) = r_2(a) = e \\
\end{split} 
\end{equation*}
One can verify that Lemmas \ref{tir}, \ref{tnr}, and \ref{tdr} are satisfied, and that $X$ is a well defined $\TL_n$-space.
\end{example}

\subsection{Long-distance intersections of the same cardinality are isomorphic} By Lemma \ref{nifi}, every intersection is either the full intersection or is a long-distance intersection, i.e is of the form $A_{i_1} \intersection ... \intersection A_{i_d}$ where $|i_k - i_l| \ge 2$ for any $k \not= l$. The natural next step is to enquire about the structural properties of the collection of long-distance intersections. The Lemma that follows shows us that the homeomorphism class of a long-distance intersection only depends on the number of subspaces intersected.

\

\begin{lemma}[Intersections of the same cardinality are homeomorphic] \label{isch} 
Let $i_1,..,i_m$ and $j_1, ... , j_m$ be such that $|i_k - i_l | \ge 2$, $|j_k - j_l | \ge 2$ for all $k \not= l$. Then,

$$A_{i_1} \intersection ... \intersection A_{i_m} \isom A_{j_1} \intersection ...\intersection A_{j_m}$$

That is, the homeomorphism class of a long-distance intersection only depends on the number of spaces that intersect. Moreover, ordering indices so that $i_1 < ... < i_m$, this isomorphism from $A_1 \intersection A_3 \intersection A_5 ... \intersection A_{2m-1}$ to $A_{i_1} \intersection .. \intersection A_{i_m}$ is explicitly given by the following map:

$$r_{i_1} r_{i_1 - 1} ... r_1r_{i_2} r_{i_2 - 1} ... r_3 ... r_{i_m} r_{i_m - 1} ... r_{2m-1} \restriction_{A_1 \intersection A_3 ... \intersection A_{2m-1}}: A_1 \intersection A_3 ... \intersection A_{2m-1} \isom A_{i_1} \intersection .. \intersection A_{i_m}$$
\end{lemma} 

\begin{proof}
Start with $A_1 \intersection A_3 \intersection ... \intersection A_{2m-1}$. We will show that any $A_{i_1} \intersection ... \intersection A_{i_m}$ is homeomorphic to $A_1 \intersection A_3 \intersection ... \intersection A_{2m-1}$. We order indices so that $i_1 < ... < i_m$. Then, observe that $i_m \ge 2m-1$. Consider $r_{i_m} r_{i_m - 1} ... r_{2m-1}$. We claim that $r_{i_m} r_{i_m - 1} ... r_{2m-1} \restriction{A_1 \intersection A_3 \intersection ... \intersection A_{2m-1}}$ is a homeomorphism from $A_1 \intersection A_3 \intersection ... \intersection A_{2m-1}$ to $A_1 \intersection A_3 \intersection ... \intersection A_{2m-3} \intersection A_{i_m}$. 

\

First, observe that the range is indeed correct:  Similar to Observation \ref{compcommrretr},  the image of $r_{i_m} r_{i_m - 1} ... r_{2m-1}  \restriction_{A_1 \intersection ... \intersection A_{2m-1}}$ is contained in each $A_{2k-1}$ for $k < m$, since, for $x \in A_1 \intersection ... \intersection A_{2m-1}$, $r_{i_m} r_{i_m - 1} ... r_{2m-1}(x) = r_{i_m} r_{i_m - 1} ... r_{2m-1} r_{2k-1}(x) = r_{2k-1}  r_{i_m} r_{i_m - 1} ... r_{2m-1} (x) \in A_{2k-1}$. Moreover, the image is contained in $A_{i_m}$ since the image of $r_{i_m}$ is $A_{i_m}$. In order to see that $A_1 \intersection A_3 \intersection ... \intersection A_{2m-3} \intersection A_{i_m}$ is the full image and that the map is a homeomorphism onto the image, we may write down the inverse, which is $r_{2m-1} r_{2m} ... r_{i_m} \restriction_{A_1 \intersection A_3 \intersection ... \intersection A_{2m-3} \intersection A_{i_m}}$. The reason why this is an inverse is that:

\begin{equation*}
\begin{split}
r_{2m-1} r_{2m} ... r_{i_m} r_{i_m} r_{i_m - 1} ... r_{2m-1} &= r_{2m-1} r_{2m} ... r_{i_m}^2 r_{i_m - 1} ... r_{2m-1} \\
&=  r_{2m-1} r_{2m} ... r_{i_m - 1} r_{i_m} r_{i_m - 1} ... r_{2m-1} \\
&=  r_{2m-1} r_{2m} ... r_{i_m - 2} r_{i_m - 1} r_{i_m -2} ... r_{2m-1} \\
&= ... \\
&= r_{2m-1} r_{2m} r_{2m-1} \\
&= r_{2m-1} \\
\end{split} 
\end{equation*}

and $r_{2m-1} \restriction A_{2m-1} = \id_{A_{2m-1}}$.

\

Therefore, we have shown that 

 $$A_1 \intersection A_3 \intersection ... \intersection A_{2m-1} \isom A_1 \intersection A_3 \intersection ... \intersection A_{2m-3} \intersection A_{i_m}$$
 
 \
 
 Repeating this process by moving $A_{2k-1}$ to $i_{k}$ for $k < m-1, m-2$ etc.. we deduce that 
 
$$A_1 \intersection A_3 \intersection ... \intersection A_{2m-1} \isom A_{i_1} \intersection ... \intersection A_{i_m}$$
\end{proof}

\

\subsection{Surjective actions and the filtration of retracts} \label{ssfr}

In the remainder of the paper, we will always work in the situation where our $\TL_n$-space $X$ is the union $X = \Union_{i=1}^{n-1} A_i$. We can think of this as ignoring ``superfluous parts of our space upon which there is no action''. 

\begin{definition}[Surjective action] \label{sa}
Let $X$ be a space on which $\TL_n$ acts. Each $u_i$ acts via a retraction onto some subspace $A_i$. We say that the action is \underline{surjective} or that $X$ is a \underline{surjective} $\TL_n$-space if $$X = \Union_{i=1}^{n-1} A_i$$
Equivalently, letting $\overset{\sim}{\TL_n} = \TL_n - \{1\}$ denote the ideal of $\TL_n$ containing all elements having cups, a surjective action is one, such that, for any $y \in X$, there exists some $x \in X$ and $a \in \overset{\sim}{\TL_n}$ such that $a \cdot x = y$, hence the name surjective.
\end{definition}

Lemma \ref{nifi} told us that any neighborhood intersection is the full intersection. Lemma \ref{isch} told us that no particular long-distance intersection is special - the intersections are all homeomorphic to intersections of the form $A_1 \intersection A_3 \intersection ... \intersection A_{2m-1}$. This motivates us to give importance to the following definition:

\

\begin{definition} (Filtration of retracts) \label{fr}
Let $X = \Union_{i=1}^{n-1} A_i$ be a surjective $\TL_n$-space (each $u_i$ acts via a retraction onto $A_i$). We define the \underline{filtration of retracts} associated to $X$, denoted by $\scrF(X)$ to be the filtration:

\begin{equation*}
\scrF(X) = A_1 \supseteq A_1 \intersection A_3 \supseteq ... \supseteq A_1 \intersection A_3 \intersection ... A_{2 \cdot \floor{\frac{n}{2}}-1} 
\end{equation*}

We will often denote $A_1 \intersection A_3 \intersection ... \intersection A_{2p-1}$ by $R_p$.
\end{definition}

\

Each $R_p$ in the above definition is a retract, and the retraction map is given by $r_1 r_3 r_5 .. r_{2p-1}$ (this is Lemma \ref{tdr}). 

\begin{remark}
The above filtration $\scrF$, defined in Definition \ref{fr} is an important definition, for two reasons: \\
(1) Theorem \ref{main} will show that the representation structure of homology groups of a $\TL_n$ space will only depend on the filtration $\scrF$, under assumptions of trivial invariants. \\
(2) $\scrF$ will provide us a natural way to define topological stability (Definition \ref{fs}). \\
\end{remark}

\section{Representation stability for Temperley-Lieb algebras at $\delta = 1$} \label{rstlad1}

In this Section, we define a notion of representation stability for Temperley-Lieb algebras, and prove that sequences of standard representations are representation stable. It is conceievable that a competing notion could be obtained by considering representations of $\TL(\infty)$ (representations of $\TL(\infty)$ are studied, for instance, in \cite{tl3}(Moore)). The notion of representation stability we will introduce will make no reference to $\TL(\infty)$, but rather will be analogous to the definition of $\FI$-modules as defined in \cite{rs1} (Church, Ellenberg, Farb). We will call our notion of a "stable module" as an "$\LS$-module", analogous to the naming convention for "$\FI$-modules". In Subsection \ref{sslstli}, we will explain why some infinite-link-state representations of $\TL(\infty)$ are not stable. 

\

\subsection{Defining representation stability} \label{ssdrs}
The definitions we provide assume that $\delta = 1$, and this is perhaps not a defect, but rather a feature of representation stability, at least from the viewpoint of actions on finite sets. With regard to topological actions, we are only interested in the $\delta = 1$ case and so we face no problems in this regard. One advantage of the notion of representation stability we provide is that it is naturally analogous to the definition of $\FI$-modules, and hence our theorem that topological stability implies representation stability (Theorem \ref{fsirs}) can be viewed as an analogue to the statement that the homology of configuration spaces is a finitely generated $\FI$-module, as in \cite{rs1} (Church, Ellenberg, Farb).

\

Let us think about the representation stability of symmetric groups as a motivation: there, one has natural sets upon which each $S_n$ acts (namely $\{1,2,...,n\}$). One then considers the chain of these sets with the natural inclusions, and considers functors from the corresponding category to the category of vector spaces (or $\Z$-modules). We will do the analogous thing here. The natural set upon which $\TL_n$ acts is the set of all link states $M_n$ (we follow the notation of \cite{tl1} (Ridout, Saint-Aubin), Section 3). That is,

$$M_n = \Union_{p=1}^{\floor{\frac{n}{2}}} \{(n,p)-\textrm{link states}\}$$

\

Adding loose strands to the end of each link state gives inclusions $M_m \inj M_n$ for $m < n$, and thus we form the chain:

$$M_1 \inj M_2 \inj M_3 ...$$

which we might compare (for analogy) with the $\FI$-module chain

$$\{1\} \inj \{1,2\} \inj \{1,2,3\} \inj...$$

Before we can define the analogous notion of an $\FI$-module, we first we define the following intermediate category  $\overset{\sim}{\LS}$:

\begin{definition}[The category $\overset{\sim}{\LS}$] \label{lstilde}
We define the category $\overset{\sim}{\LS}$ to be the category such that:

\

\underline{Objects:} The objects of  $\overset{\sim}{\LS}$ are natural numbers $n \in\N$. To emphasize they are objects, we will denote the object $n$ by $[n]$. \\

\

\underline{Morphisms:} The morphisms in $\overset{\sim}{\LS}$ are freely generated by:

\begin{itemize}
\item For any $n \in \N$, morphisms $[n] \to [n]$ are Temperley Lieb diagrams in $\TL_n$. 
\item Furthermore, for $m,n \in \N$, $m < n$, we have a single morphism $i_{m,n}: [m] \to [n]$.
\end{itemize}
\end{definition}

\begin{definition}[The category $\LS$] \label{lscat}
There is a functor from the category $F: \overset{\sim}{\LS}  \to \textbf{\emph{Set}}$, which takes each object $[n] \in \N$ to $F([n]) = M_n$, and takes, for $m < n$, each $i_{m,n}: [m] \to [n]$ to the morphism $F(i_{m,n}): M_m \inj M_n$ which adds $n-m$ loose strands to the end of each link state of $M_m$. For example, 
$$F(i_{5,7}) \left(\begin{tikzpicture}[scale = 0.2, baseline={(0,0.1)}]
    \draw (0,0) [upup=1];
    \draw (2,0) -- (2,-1);
    \draw[fill = black] (2,0) circle(0.15); 
    \draw (3,0) [upup=1];
 \end{tikzpicture}\right) =  \left(\begin{tikzpicture}[scale = 0.2, baseline={(0,0.1)}]
    \draw (0,0) [upup=1];
    \draw (3,0) [upup=1];
    \draw (2,0) -- (2,-1);
    \draw[fill = black] (2,0) circle(0.15); 
    \draw (5,0) -- (5,-1);
    \draw (6,0) -- (6,-1);
    \draw[fill = black] (5,0) circle(0.15); 
    \draw[fill = black] (6,0) circle(0.15); 
 \end{tikzpicture}\right)$$
We then define $$\LS = Image(F)$$
\end{definition}

\begin{remark}
The functor $F$ in the definition of $\LS$ above is not faithful. Different elements of $\TL_n$ can give rise to identical morphisms; for instance $\left(\begin{tikzpicture}[scale = 0.2, baseline={(0,0.1)}]
    \draw (0,2) [upup=1];
    \draw (2,2) [upup=1];
    \draw (0,0) [dndn=3];
    \draw (1,0) [sdndn=1];
  \end{tikzpicture}\right)$ and $\left(\begin{tikzpicture}[scale = 0.2, baseline={(0,0.1)}]
    \draw (0,2) [upup=1];
    \draw (2,2) [upup=1];
    \draw (0,0) [dndn=1];
    \draw (2,0) [dndn=1];
  \end{tikzpicture}\right)$ give rise to the same morphism: $M_4 \to M_4$. \\ 
\end{remark}

\begin{definition}[$\LS$-modules] \label{lsm}
An \underline{$\LS$-module} is a functor from $\LS$ to the category of vector spaces (over, say, a field $\k$).
\end{definition}

\

To elucidate the above definition, note that an $\LS$-module comprises of a chain $\{V_n\}_{n \in \N}$, where each $V_n$ is a (left) representation (or module) of $\TL_n$, and for every inclusion $F(i_{m,n}): M_m \inj M_n$, there is a $\TL$-equivariant map $(i_{m,n})_*: V_m \to V_n$ (we dropped the $F$ in the notation of $(i_{m,n})_*$ for notational ease). Of course, there are more properties which an $\LS$-module will need to satisfy, and we will discuss that in future Subsections.

\

We will end this Subsection by defining finite generation: the key features of representation stability are obtained only when the stable modules are finitely generated. For example, in the world of $\FI$-modules, the eventually-polynomial behavior of characters requires finite generation to hold. 

\

\begin{definition}[Finite generation of $\LS$-modules] \label{fglsm}
Suppose that $\{V_n\}_{n \in \N}$ is an $\LS$-module. For any collection $x_1,..,x_k \in V_{m_1},...,V_{m_k}$, we let $\spn(x_1,..,x_k)$ denote the $\LS$-submodule generated via $\TL_n$ actions on each $V_n$ and inclusions $(i_{m,n})_*$ for each $m < n$. We define the rank:

$$\rk_{\LS} (\{V_n\}_{n \in \N}) = \min\{k \sep \textrm{ there exists } x_1,..,x_k \textrm{ such that } \spn(x_1,..,x_k) =  \{V_n\}_{n \in \N} \}$$

We say that $\{V_n\}_{n \in \N}$ is finitely generated if $\rk_{\LS} (\{V_n\}_{n \in \N}) < \infty$. 
\end{definition}

\

\subsection{General form of a stable module criterion}
It is useful in representation stability to have a criterion for determining whether a chain of representations is stable or not. For example, in the case of symmetric groups $S_n$, the $\FI$-module criterion says that (see for instance Exercise 9 of \cite{rs4} (Wilson)): 

\

\begin{proposition}[The $\FI$-module criterion]
A sequence $\{V_n\}_{n \in \N}$ of symmetric group representations together with maps $\{(i_{m,n})_*: V_m \to V_n\}_{m < n}$ is an $\FI$-module if and only if 

\begin{enumerate}
\item (Compatibility of inclusions) For any $k < m < n$, $$(i_{m,n})_* \circ (i_{k,m})_* = (i_{k,n})_*$$ 
\item (Equivariance of $(i_{m,n})_*$) For any $\tau \in S_m \inj S_n$, $$\tau \circ (i_{m,n})_* = (i_{m,n})_* \circ \tau$$
\item (The $\FI$-module criterion) For any $\sigma \in S_{n-m}$, $$\sigma \o (i_{m,n})_* = (i_{m,n})_*$$
\end{enumerate}
\end{proposition}

\

We now state the general form of a stable module criterion, which will be useful for us. 

\

\begin{proposition}[General form of a stable module criterion] \label{gfsmc}
Let $\{\scrA_n\}_{n \in \N}$ be a family of monoids which include into one another with fixed inclusions $\{A_m \inj A_n \}_{m<n, m,n \in \N}$, and let $\{M_n\}_{n \in \N}$ be a collection of finite sets such that $M_n$ carries an action of $\scrA_n$, and suppose that there are $\scrA_m$-equivariant maps $i_{m,n}: M_m \inj M_n$ for $m < n$. Consider the category whose objects are $\{M_n\}_{n \in \N}$ and whose morphisms are generated by elements of $\{\scrA_n\}_{n \in \N}$ together with $\{i_{m,n}\}_{m < n}$. Suppose we define a stable module to mean a functor from this category to the category of vector spaces, and let each $M_n$ get mapped to $V_n$. Then, $\{V_n\}_{n \in \N}$ together with maps $(i_{m,n})_*: V_m \to V_n$ (the image of $i_{m,n}$) is a stable module if and only if: 

\begin{enumerate}
\item (Compatibility of inclusions) For any $k < m < n$, $$(i_{m,n})_* \circ (i_{k,m})_* = (i_{k,n})_*$$ 
\item (Equivariance of $(i_{m,n})_*$) For any $\tau \in \scrA_m \inj \scrA_n$, $$\tau \circ (i_{m,n})_* = (i_{m,n})_* \circ \tau$$
\item (The stable-module criterion) If $a, b \in \scrA_n$ are such that $a \circ i_{m,n} = b \circ i_{m,n}$, then
$$a \circ (i_{m,n})_* = b \circ (i_{m,n})_*$$
\end{enumerate}
\end{proposition}

\begin{proof}
Note that each element of $\scrA_n$ gives rise to a morphism $V_n \to V_n$, which is the image of the morphism $M_n \to M_n$ due to the action on $M_n$ (these morphisms need not be unique). We need to show that any relation between the morphisms (generated by $\{\scrA_n\}_{n \in \N}$ and $\{i_{m,n}\}_{m < n}$) can be reduced to relations of the above form.

\

Note first that (1) and (2) must hold since inclusions on $\{M_n\}_{n \in \N}$ are equivariant and compatible. Our task is therefore to show that the only other possible relation is (3).

\

By (1), and since $a,b \in \scrA_k \implies ab \in \scrA_k$, the relation must be of the form:

$$a_n i_{m_d,n} ... i_{m_2, m_3} a_{m_2} i_{m_1, m_2} a_{m_1} = a'_{n'} i_{m_d',n} ... i_{m'_2, m'_3} a'_{m'_2} i_{m'_1, m'_2} a'_{m'_1}$$

\

By (2), we may move all the $a_i$s to the left, and applying (1) again, we see that any remaining relation is of the form:

\

$$a i_{m,n} = b i_{m',n'}$$

Comparing domain and range, we see that $m = m'$ and $n = n'$. Therefore, any relation other than (1) and (2) must be of the form:

$$a i_{m,n} = b i_{m,n}$$
\end{proof}

\

\subsection{A chain of standard representations of $\TL_n$ is representation stable} The general form of a stable module criterion allows us to observe some general principles to construct stable modules. This will allow us to see that the chains of standard representations are $\LS$-modules.

\
 
\begin{proposition}[General principles for constructing stable modules]\label{gpcsm}
Let $\{\scrA_n\}_{n \in \N}$ be a family of monoids. Suppose the sequence $\{V_n\}_{n \in \N}$ together with maps $(i_{m,n})_*: V_m \to V_n$ is a stable sequence of $\scrA_n$ representations. Let $\{K_n\}_{n \in \N}$ be a sequence of subrepresentations ($K_n \le V_n$) such that $(i_{m,n})_* K_m \le K_{n}$. Then: 
\begin{enumerate}
\item $\{K_n\}_{n \in \N}$ is a stable sequence of representations. \\ 
\item $\{\quotient{V_n}{K_n} \}_{n \in \N}$ is a stable sequence of representations. \\
\end{enumerate}
\end{proposition}

\begin{proof}
 For both (1) and (2), inclusions are well defined because  $(i_{m,n})_* K_m \le K_{n}$. Inclusions are obviously compatible in both cases. Moreover, since $(i_{m,n})_*$ is equivariant as a map from $V_m$ to $V_n$, we see that this equivariance is preserved for submodules and quotient modules.
  
 \
 
 Finally, if $a \o (i_{m,n})_* = b \o (i_{m,n})_*$ holds for $\{V_n\}_{n \in \N}$, we must also have that $a \o (i_{m,n})_* \restriction_{K_m} = b \o (i_{m,n})_* \restriction_{K_m}$, and thus (1) satisfies the stable module criterion. Similarly, (2) satisfies the stable module criterion because an equality which holds in a vector space must also hold in any quotient (an equality in a quotient is a weaker statement than equality).
\end{proof}

\

We next introduce a simple shorthand notation which will allow us to freely start at any $N \in \N$ rather than at $1$, by which we mean:

\

\begin{observation}[Starting from $N$]\label{sfn}
Suppose that we are given a functor from the subcategory of $\LS$ whose objects are $\{M_n\}_{n \ge N}$ to the category of vector spaces, and let $V_n$ be the image of $M_n$. Set $V_0 = V_1 = ... = V_{N-1} = 0$, and let the images of inclusions to and from $V_k$ for $k \le N$ to be zero maps. Then, $\{V_n\}_{n \in \N}$ is a stable sequence of representations.
\end{observation}

\begin{proof}
This is true since inclusions of $0$ modules are trivially compatible, trivially and the stable module criterion trivially holds - hence, all three requirements of Proposition \ref{gfsmc} are satisfied. 
\end{proof}

\

In light of the above observation, it will be convenient for us to refer to functors from the subcategory of $\LS$ whose objects are $\{M_n\}_{n \ge N}$ to the category of vector spaces also as $\LS$-modules, and we will denote them by $\{V_n\}_{n \ge N}$.

\

For the representation stability of Temperley-Lieb algebras, the generalities above actually provide us useful information, because the representation theory of $\TL_n$ revolves around the standard representations, which are quotients of the $\C$-span of $M_n$ (which we will denote by $\C M_n$). Following the notation of  \cite{tl1} (Ridout, Saint-Aubin), we will denote the standard representation generated by $(n,p)$ link states by $V_{n,p}$. 

\

\begin{corollary}[Chains of standard representations are stable] \label{csrs}
For $p$ fixed, consider the sequence $\{V_{n,p} \}_{n \ge N}$, a sequence of standard representations of $(n,p)$ link states together with the inclusions induced by inclusions on $\{M_n\}_{n \in \N}$. Then,  $\{V_{n,p} \}_{n \ge N}$ is an $\LS$-module.
\end{corollary}

\begin{proof}
For each $p$, $V_{n,p}$ is a quotient module:

$$V_{n,p} = \quotient{\k \left(\Union_{q \ge p} \{(n,q)-\textrm{link states}\}\right)}{\k \left(\Union_{q > p} \{(n,q)-\textrm{link states}\}\right)}$$

Therefore, $V_{n,p}$ is a quotient of a submodule of $\C M_n$, and the hypotheses of Proposition \ref{gpcsm} are satisfied for both the submodule and the sub-submodule in question. Therefore, we conclude by Proposition \ref{gpcsm} that  $\{V_{n,p} \}_{n \ge N}$ is an $\LS$-module.
\end{proof}

\

Moreover, chains of standard representations have $\LS$-rank 1 (in particular, they are finitely generated).

\

\begin{proposition}[Chain of standard representations have $\LS$-rank 1] \label{csrlsr1}
$$\rk_{\LS} \{V_{n,p} \}_{n \ge N} = 1$$
\end{proposition}

\begin{proof}
Take any $(N,p)$ link state $v$ in $V_{N,p}$. For each $n \ge N$, $(i_{N,n})_* v$ is then a $(n,p)$ link state of $V_{n,p}$. Every diagram element of $V_{n,p}$ which is not in the kernel of the usual bilinear form on $V_{n,p}$ is a cyclic generator of $V_{n,p}$ (see for example \cite{tl1} (Ridout, Saint-Aubin), the proof of Proposition 3.3), and thus, when $\delta \not=0$, every link state is a cyclic generator of $V_{n,p}$ - this is because for $\delta \not=0$, the inner product of a link state with itself is some power of $\delta$ (and we are considering the case when $\delta = 1$). Therefore, by cyclicity we deduce that $\spn(v) = \{V_{n,p}\}_{n \ge N}$. Hence, $\rk_{\LS} \{V_{n,p} \}_{n \ge N} = 1$. 
\end{proof}

\

\subsection{Example: an infinite-link-state representation of $\TL(\infty)$ which is not stable} \label{sslstli}

\

Recall that, in line with Moore \cite{tl3}, $\TL(\infty)$ is defined as the the direct limit $\lim \TL_n$ with respect to the natural inclusions. As a matter of curiosity and/or completeness, a reader who is interested in stability might ask for an example of a representation of $\TL(\infty)$ which is not an $\LS$-module. We will produce such an example in this Subsection. First, we need to learn how to translate between the notion of $\TL(\infty)$-representations and the notion of $\LS$-modules. We may do that as follows: 

\

\begin{observation}[] \label{}
Suppose that $\{V_n\}_{n \in \N}$ is an $\LS$-module. Then the direct limit $\lim_n V_n$ is a representation of $\TL(\infty)$.
\end{observation}

\begin{proof}
We have $\lim_n V_n = \quotient{\ds_{n \in \N} V_n}{\sim}$ where $\sim$ is the equivelence relation $a \sim (i_{m,n})_* a$, for the maps $(i_{m,n})_*: V_m \to V_n$.  Given $a \in \TL(\infty)$ and $v \in \lim_n V_n$, define $a \cdot v$ as follows: we may choose $N$ large enough such that $a \in \TL_N \le \TL(\infty)$ and $v \in V_N$. Then define $a \cdot v$ to be via the action of $\TL_N$ on $V_N$. We must show that this is well defined (i.e does not depend on the choice of $N$). Suppose that $N_1 \le N_2$ and consider $a \in \TL_{N_1},\TL_{N_2}$ and $v \in V_{N_1}, i_{N_1,N_2} v \in V_{N_2}$. Then, by the equivariance of $(i_{N_1,N_2})_*$ (see Proposition, \ref{gfsmc}) $a (i_{N_1, N_2})_* v = (i_{N_1,N_2})_* av \sim av$. This proves well definedness of the above action.
\end{proof}

\

\begin{definition}[] \label{}
Given a representation $V$ of $\TL(\infty)$, we say that \underline{$V$ corresponds to an $\LS$-module} if there is an $\LS$-module $\{V_n\}_{n \in \N}$ such that: $V = \Union_{n \in \N} V_n$, as in the statement of the observation above.
\end{definition}

\

In \cite{tl3} (Moore), the author constructs representations of $\TL(\infty)$ which are generated by infinite link states (one is allowed infinitely many cups if one desires). He denotes the link state representation generated by a link state $w$ by $\chi(w)$ - we will instead denote it by $V(w)$. 

\

\begin{example}
(An infinite-link-state representation of $\TL(\infty)$ which is not stable). Consider the infinite link state $w = \left(\begin{tikzpicture}[scale = 0.2, baseline={(0,0.1)}]
    \draw (0,0) [upup=1];
    \draw (4,0) [upup=1];
    \draw (8,0) [upup=1];
    \draw (2,0) -- (2,-1);
    \draw[fill = black] (2,0) circle(0.15);
    \draw (3,0) -- (3,-1);
    \draw[fill = black] (3,0) circle(0.15);
    \draw (7,0) -- (7,-1);
    \draw[fill = black] (7,0) circle(0.15);
    \draw (6,0) -- (6,-1);
    \draw[fill = black] (6,0) circle(0.15);
    \draw (10,0) -- (10,-1);
    \draw[fill = black] (11,0) circle(0.15);
    \draw (11,0) -- (11,-1);
    \draw[fill = black] (10,0) circle(0.15);
    \draw[fill = black] (12.5,0) circle(0.05);
    \draw[fill = black] (13,0) circle(0.05);
    \draw[fill = black] (13.5,0) circle(0.05);
  \end{tikzpicture}\right)$. Suppose for contradiction that $V(w)$ corresponds to an $\LS$-module. Then, there is an $\LS$-module $\{V_n\}_{n \in \N}$ such that  $V(w) = \Union_{n \in \N} V_n$. So, there is some $m \in \N$ such that $w \in V_m$. But then, there is some $m' > m$ such that on the strands $m', m'+1, m'+2, m'+3$, $w$ has the pattern $\left(\begin{tikzpicture}[scale = 0.2, baseline={(0,0.1)}]
    \draw (0,0) [upup=1];
    \draw (3,0) -- (3,-1);
    \draw[fill = black] (3,0) circle(0.15);
    \draw (2,0) -- (2,-1);
    \draw[fill = black] (2,0) circle(0.15);
  \end{tikzpicture}\right)$. Let $a \in \TL_{m'+3}$ comprise of identity-map loose-strands at all positions less than $m'$, and, from $m'$ to $m'+3$ looks like:  $\left(\begin{tikzpicture}[scale = 0.2, baseline={(0,0.1)}]
    \draw (2,2) [upup=1];
    \draw (2,0) [dndn=1];
    \draw (0,0) -- (0,2);
    \draw[fill = black] (0,0) circle(0.15);
    \draw[fill = black] (0,2) circle(0.15);
    \draw (1,0) -- (1,2);
    \draw[fill = black] (1,0) circle(0.15);
    \draw[fill = black] (1,2) circle(0.15);
  \end{tikzpicture}\right)$. Let $b \in \TL_{m'+3}$ comprise of identity-map loose-strands at all positions less than $m'$, and, from $m'$ to $m'+3$ looks like: $\left(\begin{tikzpicture}[scale = 0.2, baseline={(0,0.1)}]
    \draw (2,2) [upup=1];
    \draw (0,0) [dndn=1];
    \draw (2,0) [dnup=-2];
    \draw (3,0) [dnup=-2];
  \end{tikzpicture}\right)$. Observe that $a \circ i_{m, m'+3} = b \circ i_{m,m'+3}$ but $a \circ (i_{m, m'+3})_*(w) \not= b \circ (i_{m, m'+3})_*(w) $, since $a \circ (i_{m, m'+3})_*(w)$ has $\left(\begin{tikzpicture}[scale = 0.2, baseline={(0,0.1)}]
    \draw (0,0) [upup=1];
    \draw (2,0) [upup=1];
  \end{tikzpicture}\right)$ at positions $m', m'+1, m'+2, m'+3$ whereas $b \circ (i_{m, m'+3})_*(w)$ has $\left(\begin{tikzpicture}[scale = 0.2, baseline={(0,0.1)}]
    \draw (0,0) [upup=1];
    \draw (2,0) -- (2,-1);
    \draw[fill = black] (2,0) circle(0.15);
    \draw (3,0) -- (3,-1);
    \draw[fill = black] (3,0) circle(0.15);
  \end{tikzpicture}\right)$ at positions $m', m'+1, m'+2, m'+3$. Therefore, $\{V_n\}_{n \in \N}$ does NOT satisfy the stable module criterion of \ref{gfsmc}, and hence is not an $\LS$-module. 
\end{example}

\

\section{Homology groups as $\TL_n$-representations} \label{shgtlr}

\subsection{Topological Observations} \label{sstl}

The goal of this Subsection is to make topological observations which are required to prove Theorem \ref{main}. We begin with a couple of simple but important observations. 

\begin{observation}\label{inclinjhom}
Let $r: X \surj A$ be a retraction. Then, the inclusion map $\iota: A \inj X$ induces an injection on homology $\iota_*: H(A) \inj H(X)$.
\end{observation}

\begin{proof}
$$\iota_*(\alpha) = \iota_*(\beta) \implies (r_* \iota_*) \alpha = (r_* \iota_*) \beta \implies (r \circ \iota)_* \alpha = (r \circ \iota)_* \beta \implies \alpha = \beta$$
where the last line is because $r$ is a retraction and so $r \circ \iota = \textrm{id}_A$.
\end{proof}

\

\begin{observation}[Homology classes in commuting retractions must live in their intersection] \label{htdr}
Let $X$ be a topological space, let $r: X \to A$, $r': X \to A'$ be commuting retractions, and $\tilde{\alpha} \in H_k(X)$.  Let $\iota_A, \iota_{A'}, \iota_{A \intersection A'}$ denote the inclusions $A \inj X, A' \inj X, A \intersection A' \inj X$ respectively. Suppose that there existed $\alpha \in H_k(A)$, $\alpha' \in H_k(A')$ such that $(\iota_A)_* \alpha = (\iota_{A'})_* \alpha' = \tilde{\alpha} \in H_k(X)$. Then, there is some $\beta \in H_k(A \intersection A')$ such that $(\iota_{A \intersection A'})_* \beta = \tilde{\alpha}$.
\end{observation}

\begin{proof}
Since $r$ and $r'$ commute, $r' r$ is a retraction onto $A \intersection A'$ (Observation \ref{compcommrretr}). Therefore, observe that $(r r')_* \tilde{\alpha} \in H_k(A \intersection A')$. Set $\beta = (r r')_* \tilde{\alpha}$.

\

Observe that $(\iota_A \circ r) \circ (\iota_{A'} \circ r') = \iota_{A \intersection A'} \circ (rr')$ as maps from $X$ to $X$. Hence, $(\iota_{A \intersection A'})_* \beta = (\iota_{A \intersection A'} \circ r r')_* \tilde{\alpha} = (\iota_A \circ r \circ \iota_{A'} \circ r')_*  (\iota_{A'})_* \alpha' = (\iota_A \circ r)_* \tilde{\alpha} = (\iota_A \circ r)_* (\iota_A)_* \alpha = \tilde{\alpha}$. And so we are done.
\end{proof}

\

We now define the notion of a minimal intersection. This notion as has a twofold importance: \\
(1) We will use it in Subsection \ref{sscmgpi} to motivate the map that identifies the cyclic module generated by a homology class with the standard representation. \\
(2) It is the conceptual reason for the injectivity of that map. \\ 

\begin{definition}
Let $X$ be a surjective $\TL_n$-space. Let $H_k(X)^{\TL_n}$ denote the invariant subrepresentation of $H_k(X)$, that is, $H_k(X)^{\TL_n} = \{\alpha \in H_k(X) \sep u_i \alpha = \alpha \textrm{ } \forall i\}$.
\end{definition}

\begin{definition}(Minimal Intersection) \label{defminimalintersection}
Let $X$ be a surjective $\TL_n$-space. Given an element $\alpha \in H_k(X)$, and given  $\{i_1,..,,i_m\}$ such that $|i_j - i_l| \ge 2$ for all $j \not= l$, we say that $A_{i_1} \intersection ... \intersection A_{i_m}$ is a \\ \underline{minimal intersection containing $\alpha$} if: \\
(1) $\alpha$ lies in the image $H_k(A_{i_1} \intersection ... \intersection A_{i_m}) \inj H_k(X)$ \\
(2) Given any $l \not\in \{i_1,...,i_m\}$, $\alpha$ does not lie in the image $H_k(A_{i_1} \intersection ... \intersection A_{i_m} \intersection A_l) \inj H_k(X)$ \\
(3) $\alpha \not\in H_k(X)^{\TL_n} \subseteq H_k(X)$
\end{definition}

The utility of minimal intersections is that they have a uniqueness property, which we will see by exploiting Observation \ref{htdr} and Corollary \ref{hnifi}. 

\

\begin{lemma}[Uniqueness of minimal intersection containing a given homology class] \label{uniquenessminimalintersection}
Let $X$ be a surjective $\TL_n$-space. Take $\alpha \in H_k(X)$, and suppose that $A_{i_1} \intersection ... \intersection A_{i_p}$ is a minimal intersection containing $\alpha$. Then, no other minimal intersection contains $\alpha$.
\end{lemma}

\begin{proof}
Suppose for contradiction that $\alpha$ is contained in two minimal intersections  $A_{i_1} \intersection ... \intersection A_{i_p}$ and $A_{j_1} \intersection ... \intersection A_{j_q}$ with $\{i_1, ... i_p\} \not= \{j_1,...,j_q\}$. 

\

Suppose first that there exists some $m$ such that $|j_m - i_l| \ge 2$ for all $l$. Then, since the retraction $u_{j_m}$ commutes with the retractions $\{u_{i_l}\}_{l=1}^p$, Observation \ref{htdr} implies that $\alpha$ lives in the image of the map

$$H_k(A_{i_1} \intersection ... \intersection A_{i_p} \intersection A_{j_m}) \inj H_k(A_{i_1} \intersection ... \intersection A_{i_p}) \inj H_k(X)$$ 

Explicitly, the proof of Observation \ref{htdr} shows that if $\beta \in H_k(A_{i_1} \intersection ... \intersection A_{i_p})$ is an element of the preimage of $\alpha$ in $H_k(A_{i_1} \intersection ... \intersection A_{i_p})$, then $u_{j_m} \beta \in H_k(A_{i_1} \intersection ... \intersection A_{i_p} \intersection A_{j_m})$ is in the preimage of $\alpha$ in $H_k(A_{i_1} \intersection ... \intersection A_{i_p} \intersection A_{j_m})$.

\

This contradicts the fact that $A_{i_1} \intersection ... \intersection A_{i_p}$ is a minimal intersection of $\alpha$ since we have shown that $\alpha$ lies in the image $H_k(A_{i_1} \intersection ... \intersection A_{i_p} \intersection A_{j_m}) \inj H_k(X)$, and hence there cannot exist $m$ such that $|j_m - i_l| \ge 2$ for all $l$.

\

The remaining case is when, for each $j_m$, there is some $i_l$ such that $|j_m - i_l| \le 1$. Since $\{i_1, ... i_p\} \not= \{j_1,...,j_q\}$, there must exist some $j_m, i_l$ such that $|j_m - i_l| = 1$. Thus, there is a $j$ such that $\alpha = \iota_* \beta_j \in H_k(A_j) \inj H_k(X)$ and $\alpha = \iota_* \beta_{j+1} \in H_k(A_{j+1}) \inj H_k(X)$. Notice then that $u_j \cdot \alpha = u_j \cdot \iota_* \beta_j = (u_j \circ \iota)_* \beta_j = \beta_j = \alpha$, and similarly $u_{j+1} \alpha = \alpha$. Thus, by the cycle map trick (Corollary \ref{hnifi}), we see that $u_i \cdot \alpha = \alpha$ for all $i$, and so $\alpha \in H_k(X)^{\TL_n}$, a contradiction to the third minimal intersection assumption. 
\end{proof}

\

An easy reformulation of the above into a form which is useful is: 

\

\begin{corollary}\label{tlc}
If $\alpha, \beta \in H_k(X)$ are such that $A_{i_1} \intersection ... \intersection A_{i_m}$ is the minimal intersection containing $\alpha$ and $A_{j_1} \intersection ... \intersection A_{j_d}$ is the minimal intersection containing $\beta$, with $\{i_1,..,i_m\} \not= \{j_1,..,j_d\}$, then $\alpha \not= \beta$.
\end{corollary}

\begin{proof}
If $\alpha = \beta$, then we have an element which has two minimal intersections, namely $A_{i_1} \intersection ... \intersection A_{i_m}$ and $A_{j_1} \intersection ... \intersection A_{j_d}$, thereby contradicting Lemma \ref{uniquenessminimalintersection}. 
\end{proof}

\begin{remark}
A somewhat subtle point here is that it is implicit that $\alpha, \beta \not\in H_k(X)^{\TL_n}$ due to the third assumption of minimal intersection in Definition \ref{defminimalintersection}. Without this crucial property of $\alpha, \beta$, the above is false, see for example Example \ref{reasonfortrivialinvariants} where $A_1 \intersection A_4 = A_1 \intersection A_3 = Q$, and so, in the notation of that example, the element of $H_1(Q) \inj H_1(X)$ has multiple minimal intersections, for example  $A_1 \intersection A_4$ and $A_1 \intersection A_3$. 
\end{remark}

\

We end this Section with a simple proposition which shows us that under ``nice enough $\TL_n$-moves, homology classes do not change much'': 

\

\begin{proposition}[Isomorphisms preserve the null-homologous property] \label{zp}
Suppose that $\alpha \in \iota_* H_k(A_{i_1} \intersection A_{i_2} \intersection ... \intersection A_{i_q}) \subseteq H_k(X)$, where $|i_l - i_m| \ge 2$ for any $l \not= m$. Suppose that $j$ is such that $| j -i_s| \le 1$ for some $s$. Then, 

$$u_j \alpha = 0 \iff \alpha = 0$$
\end{proposition}

\begin{proof}
Note that the inclusion map provides an injection $H_k(A_{i_1} \intersection A_{i_2} \intersection ... \intersection A_{i_q}) \inj H_k(A_{i_s})$. Write $\alpha = \iota_* \beta$ for $\beta \in H_k(A_{i_s})$. Since $| j -i_s| \le 1$,  $u_j \restriction_{A_{i_s}}: A_{i_s} \to A_j$ is a homeomorphism. Therefore $0 \not= \beta \in H_k(A_{i_s}) \iff 0 \not= u_j \beta \in H_k(A_j)$. Since the inclusion maps provide injections $H_k(A_{i_s})  \inj H_k(X)$, $H_k(A_{j})  \inj H_k(X)$, we conclude that $u_j \alpha \not= 0 \iff \alpha \not= 0$.
\end{proof}

\

\subsection{The reason for the assumption of trivial invariants $(H_k(X)^{\TL_n} = 0)$} \label{ssrcxq}

\

Let $X$ be a surjective $\TL_n$-space. When $H_k(X)^{\TL_n} \not= 0$, we have much less control over the $\TL_n$-representation structure of $H_*(X)$. We will now provide an example of why this is the case:

\

\begin{example}\label{reasonfortrivialinvariants}
Let $X$ be the wedge of 5 copies of $S^1$, denoted by $Q, S_1, S_2, S_3, S_4$. Define the subspaces $A_1,A_2,A_3,A_4 \subset X$ by $A_i = Q \twedge S_i$ for each $i$. Note that $A \isom S^1 \twedge S^1$ for each $i$. Note that $\intersection_{i=1}^4 A_i = Q$.  Denote $\{*\}$ be the intersection of the circles. For each $i = 1,2,3$, fix isomorphisms $\phi_i: S_i \isom S_{i+1}$ which preserve $S_i \intersection S_{i+1} = \{*\}$.

\begin{figure}[H]
\centering
\begin{tikzpicture}
   \begin{polaraxis}[grid=none, axis lines=none]
     \addplot[mark=none, color = gray, very thick, domain=-18:18,samples=300] {-cos(5*x)};
   \end{polaraxis}
   \begin{polaraxis}[grid=none, axis lines=none]
     \addplot[mark=none, color = black, very thick, domain=18:54,samples=300] {-cos(5*x)};
   \end{polaraxis}
   \begin{polaraxis}[grid=none, axis lines=none]
     \addplot[mark=none, color = black, very thick, domain=54:90,samples=300] {-cos(5*x)};
   \end{polaraxis}
    \begin{polaraxis}[grid=none, axis lines=none]
     \addplot[mark=none, color = black, very thick, domain=90:126,samples=300] {-cos(5*x)};
   \end{polaraxis}
    \begin{polaraxis}[grid=none, axis lines=none]
     \addplot[mark=none, color = black, very thick, domain=126:162,samples=300] {-cos(5*x)};
   \end{polaraxis}
   \node at (0,3.4) {$\mathbf{Q}$};
   \node at (2.8,6.8) {$\mathbf{S_1}$};
   \node at (2.8,0) {$\mathbf{S_4}$};
   \node at (6.2,5.6) {$\mathbf{S_2}$};
   \node at (6.3,1.4) {$\mathbf{S_3}$};
\end{tikzpicture}
\caption{An illustration of X. We will describe the maps $u_1,u_2,u_3,u_4$ below, which will make $X$ into a $\TL_5$-space, with $A_i = Q \twedge S_i$ for each $i$ as described above. Note that $H_1(X)^{\TL_n} =H_1(Q) =  H_1(S^1) \not= 0$}.
\end{figure}

Let $u_1,u_2,u_3,u_4$ act as follows: $u_i$ is the identity on $A_i = Q \twedge S_i$, so $u_i \restriction_Q = \textrm{id}_Q$, $u_i \restriction_{S_i} = \textrm{id}_{S_i}$. If $i > 1$, $u_i \restriction_{S_{i-1}} = \phi_{i-1}$. If $i < 4$, $u_i \restriction_{S_{i+1}} = \phi_i^{-1}$. Fix an isomorphism $\psi_4: S_4 \isom Q$ which fixes $\{*\}$. Define $\psi_3 := \psi_4 \circ \phi_3: S_3 \isom Q$ and $\psi_1 := \psi_3 \circ \phi_2 \circ \phi_1: S_1 \isom Q$. We define the long distance maps as follows: 

\begin{equation*}
\begin{split}
& u_1(S_3) = u_1(S_4) = Q, \textrm{via } u_1 \restriction_{S_3} = \psi_3, u_1 \restriction_{S_4} = \psi_4 \\ 
& u_2(S_4) = \{*\} \\ 
& u_3(S_1) = Q, \textrm{via } u_3 \restriction_{S_1} = \psi_1 \\ 
& u_4(S_1) = u_4(S_2) = \{*\} \\ 
\end{split} 
\end{equation*} 

That is, $u_2$ and $u_4$, when applied long distance, retract to $\{*\}$ rather than the whole of $Q$. Let us observe that this is a well defined $\TL_n$-action. Neighbor relations (Lemma \ref{tnr}) hold because the maps, when restricted to neighboring copies of $S^1$, have been chosen to be mutually compatible isomorphisms. What remains to be checked are the long-distance relations. There are three cases to consider:

\

(1) To show that $u_1 u_4 = u_4 u_1$: If $x \in Q \union S_1 \union S_2$, then both maps act on $x$ by the trivial map to $\{*\}$. If $x \in S_3$, then $u_1 u_4(x) = \psi_4 \phi_3(x) = \psi_3(x) = u_1(x) = u_4 u_1(x)$, where the last equality follows because $u_1(x) \in Q$ and $u_4$ fixes $Q$. If $x \in S_4$, then $u_4 u_1(x) = u_4 \psi_4(x) = \psi_4(x)$ and $u_1 u_4(x) = u_1 \textrm{id}(x) = u_1(x)  = \psi_4(x)$. Thus, in all cases, $u_1 u_4 = u_4 u_1$.

\

(2) To show that $u_1 u_3 = u_3 u_1$:  If $x \in S_1$, $u_1 u_3 (x) = u_1 \psi_1(x) = \psi_1(x) = \psi_1 (\textrm{id}(x)) = u_3 u_1(x)$. Similarly, if $x \in S_3$, $u_1 u_3(x) = \psi_3(x) = u_3 u_1(x)$. If $x \in S_2$, then 
$u_1 u_3(x) = u_1 \phi_2(x) = \psi_3 \phi_2(x) = \psi_3 \phi_2 \phi_1 \phi_2(x) = \psi_3 \phi_2(x) = u_1 u_3(x)$. Finally, if $x \in S_4$, then $u_3 u_1(x) = u_3 \psi_4(x) = \psi_4(x) = \psi_3 \phi^{-1}_3(x) = u_1 u_3(x)$

\

(3) To show that $u_2 u_4 = u_4 u_2$: For $x \in S_1 \union S_2 \union S_3 \union S_4$, note that $u_4(x) \in S_4$ and so $u_2 u_4(x) \in \{*\}$. Likewise, $u_2(x) \in S_2$ and so $u_4 u_2(x) \in \{*\}$.

\

Therefore, the above action of $\TL_5$ on $X$ is well defined. However, the $\TL_5$-representation structure of $H_1(X)$ is not a direct sum of standard representations. In fact, the representation $H_1(X)$ is a nontrivial extension $0 \to [1] \to H_1(X) \to V_{5,1} \to 0$, where $[1]$ the one dimensional representation of $\TL_5$ wherein every element of $\TL_n$ fixes the spanning vector. To see this, note first that since $Q = X^{\TL_n}$, the fixed points of $\TL_n$ action, $\TL_n$ acts on $\quotient{X}{Q}$. Thus, consider the short exact sequence $0 \to Q \inj X \surj \quotient{X}{Q} \to 0$. Each map here is $\TL_n$-equivariant and this is therefore a short exact sequence of $\TL_n$-spaces.  This gives rise to a long exact sequence in homology $$0 \to H_0(Q) \to H_0(X) \to H_0(\quotient{X}{Q}) \to H_1(Q) \to H_1(X) \to H_1(\quotient{X}{Q}) \to H_2(Q) \to ...$$

Since  $H_0(X) \to H_0(\quotient{X}{Q})$ is a surjection, $H_0(\quotient{X}{Q}) \to H_1(Q)$ is the zero map, and $H_2(Q) = 0$. So we get the exact sequence of Homology groups as $\TL_n$-representations, given by $0 \to  H_1(Q) \to H_1(X) \to H_1(\quotient{X}{Q}) \to 0$. Note that $H_1(Q) = [1]$, $H_1(\quotient{X}{Q})  \isom V_{5,1}$ (it is a simple case of Example \ref{ex1}). Therefore, $H_1(X)$ is an extension $0 \to [1] \to H_1(X) \to V_{5,1} \to 0$. 

\

This extension does not split because, if it did, it would be left split, and so there would have to be a map $t: H_1(X) \to H_1(Q)$ such that $t \circ i = \id_Q$.  $\{Q,S_1,S_2,S_3,S_4\}$ forms a basis of $H_1(X)$. Note that $t(u_3 S_1) = t(Q) = Q$ and so by equivariance, $u_3 t(S_1) = Q$, and so $t(S_1) = Q$. On the other hand, $t(u_4 S_1) = t(0) = 0$, and so $u_4 t(S_1) = 0$ and so $t(S_1) = 0$. This is a contradiction since $0 \not= Q \in H_1(X)$. Thus, $H_1(X)$ is a nontrivial extension $0 \to [1] \to H_1(X) \to V_{5,1} \to 0$ (does not split).

\

We may consider another action of $\TL_5$ on $X$ by taking all long distance maps to be trivial maps onto $\{*\}$ away from $Q = \Intersection_{i=1}^4 A_i$ (note that $Q$ must be fixed), and this would have given us another $\TL_5$-space with the same filtration but with a different representation structure for $H_1(X)$ (the new representation structure would be given by $V_{5,1} \ds [1]$). 

\

Thus, since we constructed two different actions on $X$, with the same filtration but very different representation structure of homology, we conclude that the representation structure of $H_1(X)$ could not be read off from the filtration. Thus, if we do not require that $H_k(X)^{\TL_n} = 0$, a result like Theorem \ref{main} cannot hold as stated, and hence this is a crucial assumption to the theorem.
 \end{example}

\subsection{The cyclic module generated from a $p$-intersection} \label{sscmgpi}

In this Subsection, we will use results from Subsection \ref{sstl} to deduce a key stepping stone to Theorem \ref{main}. Before we delve into the theorem of this Subsection, we will remind the reader about a standard construction one encounters when studying the representation theory of Temperley-Lieb algebras, especially from a diagrammatic perspective: 

\

Given any element of $\TL_n$, we may construct this element by placing two link states ``back to back'' , as per the following definition:

\begin{definition}\label{vw*}
Given link states $v$ and $w$ which have the same number of cups, we define $v w^* \in \TL_n$ to be the link state where $w$ is reflected and placed below $v$. This construction is often alternatively denoted by $|v \textrm{ } w|$, as in \cite{tl1} (Ridout, Saint-Aubin)). 

\

For example 

$$\left(\begin{tikzpicture}[scale = 0.2, baseline={(0,0.1)}]
    \draw (-1,0) [upup=1];
    \draw (1,0) [upup=1];
    \draw[fill = black] (3,0) circle(0.15);
    \draw[fill = black] (4,0) circle(0.15);
    \draw (3,0) -- (3,-1);
    \draw (4,0) -- (4,-1);
  \end{tikzpicture} \right) \cdot \left(\begin{tikzpicture}[scale = 0.2, baseline={(0,0.1)}]
    \draw (1,0) [upup=3];
    \draw (2,0) [supup=1];
    \draw[fill = black] (-1,0) circle(0.15);
    \draw[fill = black] (0,0) circle(0.15);
    \draw (-1,0) -- (-1,-1);
    \draw (0,0) -- (0,-1);
  \end{tikzpicture} \right)^* = \begin{tikzpicture}[scale = 0.2, baseline={(0,0.1)}]
    \draw (-1,3) [upup=1];
    \draw (1,3) [upup=1];
    \draw[fill = black] (3,3) circle(0.15);
    \draw[fill = black] (4,3) circle(0.15);
    \draw (3,3) -- (3,2.3) -- (-1,1.4) -- (-1,0);
    \draw (4,3) -- (4,2.3) -- (0,1.4) -- (0,0);
    \draw[fill = black] (-1,0) circle(0.15);
    \draw[fill = black] (0,0) circle(0.15);
    \draw (1,0) [dndn=3];
    \draw (2,0) [sdndn=1];
  \end{tikzpicture} $$ 
  
\end{definition}

\

We now state and prove the theorem of this Subsection:

\

\begin{theorem}[Cyclic module generated from a $p$-intersection]\label{cmt}
Let $X = \union_{i=1}^{n-1} A_i$ be a surjective $\TL_n$-space such that $H_k(X)^{\textrm{TL}_n} = 0$. Consider $R_p = A_1 \intersection A_3 \intersection ... \intersection A_{2p-1}$. Suppose that: \\
(1) $0 \not= \alpha \in \iota_* H_k(R_p, \Z) \subset H_k(X, \Z)$. \\
(2) Let $\{i_1,...,i_d\}$ be a set such that  $|i_k -i_l| \ge 2$ for all $k,l$. If $\{1,3,5,...,2p-1 \} \subsetneq \{i_1,...,i_d\}$, $a \in \TL_n$ and $a \cdot \alpha \in \iota_* H_k(A_{i_1} \intersection ... \intersection A_{i_d})$, then $a \cdot \alpha = 0$. \\
(3) $m$ is the smallest natural number such that $m \alpha = 0$ (i.e $m$ records torsion, $m = \infty$ if there is no torsion). \\
Then, the cyclic submodule of $H_k(X, \Z)$ generated by $\alpha$ is the standard representation with coefficients that record the torsion:
 
 $$\TL_n \cdot \alpha = \begin{cases} V_{n,p} \x_\Z \quotient{\Z}{m \Z} & \textrm{ if } m < \infty \\ V_{n,p} & \textrm{ if } m = \infty \end{cases} $$
 
 In particular, if $\k$ is a field containing $\Q$, then, for $\alpha \in H_k(X,\k)$, $\TL_n \cdot \alpha = V_{n,p}$ is the standard representation over $\k$. 
\end{theorem}

\begin{proof}
For simplicity, we will assume that $\k$ is a field containing $\Q$ and $\alpha \in H_k(S,\k)$, since this is where all the content of the above statement is. The result over $\Z$ is easy to deduce from the statement for $\k$, and we will remark in the end why the result over $\Z$ holds.

\

We construct the map $h: V_{n,p} \to TL_n \cdot \alpha$ by: 

$$\begin{tikzpicture}[scale = 0.2, baseline={(0,0.1)}]
    \draw (0,0) [upup=1];
    \draw (2,0) [upup=1];
    \draw[fill = black] (4.5,0) circle(0.05);
    \draw[fill = black] (5,0) circle(0.05);
    \draw[fill = black] (5.5,0) circle(0.05);
    \draw (7,0) [upup=1];
    \draw (9,0) -- (9,-1);
    \draw[fill = black] (9,0) circle(0.15);
    \draw (10,0) -- (10,-1);
    \draw[fill = black] (10,0) circle(0.15);
    \draw[fill = black] (11.5,0) circle(0.05);
    \draw[fill = black] (12,0) circle(0.05);
    \draw[fill = black] (12.5,0) circle(0.05);
    \draw[fill = black] (14,0) circle(0.15);
    \draw (14,0) -- (14,-1);
  \end{tikzpicture} \mapsto \alpha$$ 

\

and we extend the map to all of $V_{n,p}$ by $\TL_n$-equivariance. We must show that this gives rise to a well defined map $h: V_{n,p} \to \TL_n \cdot \alpha$, and that this map is injective. Note that $\TL_n$-equivariance immediately implies surjectivity of the map.  

\

\underline{Motivation and connection with Subsection \ref{sstl}:} Assume for the moment that $h$ is well defined. We begin with the following observation which is the motivation for constructing the map above, and will also help us to prove injectivity later on: 

\fbox{
\begin{minipage}{\textwidth}
Let $v \in V_{n,d}$ be a nonzero link state with length one cups $\begin{tikzpicture}[scale = 0.2, baseline={(0,0.1)}]
    \draw (0,1) [upup=1];
  \end{tikzpicture}$ at positions $i_1,..,i_q$ (and no other positions). Then, the minimal intersection containing $h(v)$ is $A_{i_1} \intersection ... \intersection A_{i_q}$. 
\end{minipage}
}
  
 \
  
To see this, first observe that $h(v)$ is contained in this intersection since each $u_{i_j}$ fixes $h(v)$, and thus $h(v) =  u_{i_j} h(v) \in \iota_* H_k(A_{i_j})$. Thus, by Lemma \ref{htdr}, $h(v) = u_{i_1} u_{i_2} ... u_{i_q} h(v) \in \iota_* H_k(A_{i_1} \intersection ... \intersection A_{i_q})$. Next note that any other $u_{j}$ does not fix $v$. Note that for such a $j$, we have that $|j - i_s| > 0$ for all $s$. Suppose that we still have that $u_j h(v) = h(v)$ despite the fact that $u_j$ does not fix $v$. If $|j-i_s| = 1$ for some $s$, then by Corollary \ref{hnifi} (cycle map trick), we must have that $h(v) \in H_k(X)^{\textrm{TL}_n}$ and since we assumed that $H_k(X)^{\textrm{TL}_n} = 0$, we have that $h(v) = 0$. Since every link state is a cyclic generator of $V_{n,p}$ (since none of them lie in the kernel of the bilinear form, see for example \cite{tl1} (Ridout, Saint-Aubin), the proof of Proposition 3.3), there is some $a \in \TL_n$ such that $a v = (\begin{tikzpicture}[scale = 0.2, baseline={(0,0.1)}]
    \draw (0,0) [upup=1];
    \draw (2,0) [upup=1];
    \draw[fill = black] (4.5,0) circle(0.05);
    \draw[fill = black] (5,0) circle(0.05);
    \draw[fill = black] (5.5,0) circle(0.05);
    \draw (7,0) [upup=1];
    \draw (9,0) -- (9,-1);
    \draw[fill = black] (9,0) circle(0.15);
    \draw (10,0) -- (10,-1);
    \draw[fill = black] (10,0) circle(0.15);
    \draw[fill = black] (11.5,0) circle(0.05);
    \draw[fill = black] (12,0) circle(0.05);
    \draw[fill = black] (12.5,0) circle(0.05);
    \draw[fill = black] (14,0) circle(0.15);
    \draw (14,0) -- (14,-1);
  \end{tikzpicture})$. But then, $0 = a h(v) = h(av) = h(\begin{tikzpicture}[scale = 0.2, baseline={(0,0.1)}]
    \draw (0,0) [upup=1];
    \draw (2,0) [upup=1];
    \draw[fill = black] (4.5,0) circle(0.05);
    \draw[fill = black] (5,0) circle(0.05);
    \draw[fill = black] (5.5,0) circle(0.05);
    \draw (7,0) [upup=1];
    \draw (9,0) -- (9,-1);
    \draw[fill = black] (9,0) circle(0.15);
    \draw (10,0) -- (10,-1);
    \draw[fill = black] (10,0) circle(0.15);
    \draw[fill = black] (11.5,0) circle(0.05);
    \draw[fill = black] (12,0) circle(0.05);
    \draw[fill = black] (12.5,0) circle(0.05);
    \draw[fill = black] (14,0) circle(0.15);
    \draw (14,0) -- (14,-1);  \end{tikzpicture}) = \alpha$, a contradiction. Therefore, $|j - i_s| \ge 2$ for every $s$. If $v$ has loose strands in positions $j$ and $j+1$, $u_j v = 0$ in which case $u_j h(v) = h(u_j v) = h(0) = 0$. Yet $h(v)$ cannot be in $\iota_* H_k(A_j)$, because if it was, $h(v) = u_j h(v) = 0$, and once again, we may exploit cyclicity to see that $\alpha = 0$, a contradiction. The only remaining case is when $u_j$ ``breaks a nest''. This is best explained by means of an example, which easily generalizes: 
  
\

Suppose that $$u_7 h(\begin{tikzpicture}[scale = 0.2, baseline={(0,0.1)}]
    \draw (0,0) [bbupup=7];
    \draw (1,0) [bupup=5];
    \draw (2,0) [upup = 3];
    \draw (3,0) [supup = 1];
    \end{tikzpicture}) = h(\begin{tikzpicture}[scale = 0.2, baseline={(0,0.1)}]
    \draw (0,0) [bbupup=7];
    \draw (1,0) [bupup=5];
    \draw (2,0) [upup = 3];
    \draw (3,0) [supup = 1];
    \end{tikzpicture})$$
\

That is,

 $$h(\begin{tikzpicture}[scale = 0.2, baseline={(0,0.1)}]
    \draw (0,0) [upup=1];
    \draw (2,0) [bupup=3];
    \draw (3,0) [upup = 1];
    \draw (6,0) [upup = 1];
    \end{tikzpicture}) = h(\begin{tikzpicture}[scale = 0.2, baseline={(0,0.1)}]
    \draw (0,0) [bbupup=7];
    \draw (1,0) [bupup=5];
    \draw (2,0) [upup = 3];
    \draw (3,0) [supup = 1];
    \end{tikzpicture})$$

\

Applying $u_6$ to both sides, we have that 

 $$h(\begin{tikzpicture}[scale = 0.2, baseline={(0,0.1)}]
    \draw (0,0) [upup=1];
    \draw (2,0) [bupup=5];
    \draw (3,0) [upup = 1];
    \draw (5,0) [upup = 1];
    \end{tikzpicture}) = h(\begin{tikzpicture}[scale = 0.2, baseline={(0,0.1)}]
    \draw (0,0) [bupup=7];
    \draw (1,0) [upup=1];
    \draw (3,0) [upup = 1];
    \draw (5,0) [upup = 1];
    \end{tikzpicture})$$

\

But the left hand side is in $\iota_* H_k(A_1)$, while the right hand side is in $\iota_* H_k(A_2)$. Thus, by Corollary \ref{hnifi}, we have that $h(\begin{tikzpicture}[scale = 0.2, baseline={(0,0.1)}]
    \draw (0,0) [upup=1];
    \draw (2,0) [bupup=5];
    \draw (3,0) [upup = 1];
    \draw (5,0) [upup = 1];
    \end{tikzpicture}) \in H_k(X)^{\textrm{TL}_n}$, and thus, since $H_k(X)^{\textrm{TL}_n} = 0$, $h(\begin{tikzpicture}[scale = 0.2, baseline={(0,0.1)}]
    \draw (0,0) [upup=1];
    \draw (2,0) [bupup=5];
    \draw (3,0) [upup = 1];
    \draw (5,0) [upup = 1];
    \end{tikzpicture}) = 0$. By cyclicity, we once again have that $\alpha = 0$, a contradiction. This establishes the fact that, for a link state $v$, the minimal intersection of $h(v)$ can be read off from the position of the length $1$ cups of $v$. This motivates the construction of the above map. We now prove well-definedness: 

\

\underline{Well definedness:} We must show that if $a,b \in \TL_n$ are such that $$a\left(\begin{tikzpicture}[scale = 0.2, baseline={(0,0.1)}]
    \draw (0,0) [upup=1];
    \draw (2,0) [upup=1];
    \draw[fill = black] (4.5,0) circle(0.05);
    \draw[fill = black] (5,0) circle(0.05);
    \draw[fill = black] (5.5,0) circle(0.05);
    \draw (7,0) [upup=1];
    \draw (9,0) -- (9,-1);
    \draw[fill = black] (9,0) circle(0.15);
    \draw (10,0) -- (10,-1);
    \draw[fill = black] (10,0) circle(0.15);
    \draw[fill = black] (11.5,0) circle(0.05);
    \draw[fill = black] (12,0) circle(0.05);
    \draw[fill = black] (12.5,0) circle(0.05);
    \draw[fill = black] (14,0) circle(0.15);
    \draw (14,0) -- (14,-1);
  \end{tikzpicture} \right) = b\left( \begin{tikzpicture}[scale = 0.2, baseline={(0,0.1)}]
    \draw (0,0) [upup=1];
    \draw (2,0) [upup=1];
    \draw[fill = black] (4.5,0) circle(0.05);
    \draw[fill = black] (5,0) circle(0.05);
    \draw[fill = black] (5.5,0) circle(0.05);
    \draw (7,0) [upup=1];
    \draw (9,0) -- (9,-1);
    \draw[fill = black] (9,0) circle(0.15);
    \draw (10,0) -- (10,-1);
    \draw[fill = black] (10,0) circle(0.15);
    \draw[fill = black] (11.5,0) circle(0.05);
    \draw[fill = black] (12,0) circle(0.05);
    \draw[fill = black] (12.5,0) circle(0.05);
    \draw[fill = black] (14,0) circle(0.15);
    \draw (14,0) -- (14,-1);
  \end{tikzpicture} \right)$$
  then in fact  $a \alpha = b \alpha$. 
  
\

There are two cases to consider. The first case is when $a \left( \begin{tikzpicture}[scale = 0.2, baseline={(0,0.1)}]
    \draw (0,0) [upup=1];
    \draw (2,0) [upup=1];
    \draw[fill = black] (4.5,0) circle(0.05);
    \draw[fill = black] (5,0) circle(0.05);
    \draw[fill = black] (5.5,0) circle(0.05);
    \draw (7,0) [upup=1];
    \draw (9,0) -- (9,-1);
    \draw[fill = black] (9,0) circle(0.15);
    \draw (10,0) -- (10,-1);
    \draw[fill = black] (10,0) circle(0.15);
    \draw[fill = black] (11.5,0) circle(0.05);
    \draw[fill = black] (12,0) circle(0.05);
    \draw[fill = black] (12.5,0) circle(0.05);
    \draw[fill = black] (14,0) circle(0.15);
    \draw (14,0) -- (14,-1);
  \end{tikzpicture} \right) = 0 \in V_{n,p}$. In this case, if we think of the link state in $M_n$ (see Subsection \ref{ssdrs} if the notation is unclear), it has more than $p$ cups. Such a link state can be taken to the link state comprising of $q > p$ cups followed by $n-q$ loose strands by a sequence of neighbor retractions in the sense of Proposition \ref{zp}. Call the product of neighbor retractions $a'$. Then, in light of Observation \ref{htdr}, $a' a \alpha = 0$ by hypothesis of our Theorem, since $a' a \alpha \in H_k(R_p \intersection A_{2p+1})$. Therefore, by Proposition \ref{zp}, $a \alpha = 0$. Similarly, since $b \left( \begin{tikzpicture}[scale = 0.2, baseline={(0,0.1)}]
    \draw (0,0) [upup=1];
    \draw (2,0) [upup=1];
    \draw[fill = black] (4.5,0) circle(0.05);
    \draw[fill = black] (5,0) circle(0.05);
    \draw[fill = black] (5.5,0) circle(0.05);
    \draw (7,0) [upup=1];
    \draw (9,0) -- (9,-1);
    \draw[fill = black] (9,0) circle(0.15);
    \draw (10,0) -- (10,-1);
    \draw[fill = black] (10,0) circle(0.15);
    \draw[fill = black] (11.5,0) circle(0.05);
    \draw[fill = black] (12,0) circle(0.05);
    \draw[fill = black] (12.5,0) circle(0.05);
    \draw[fill = black] (14,0) circle(0.15);
    \draw (14,0) -- (14,-1);
  \end{tikzpicture} \right) = 0$, we have that $b \alpha = 0$. So, $a \alpha = 0 = b \alpha$.

\

\

The remaining case is when $a \left( \begin{tikzpicture}[scale = 0.2, baseline={(0,0.1)}]
    \draw (0,0) [upup=1];
    \draw (2,0) [upup=1];
    \draw[fill = black] (4.5,0) circle(0.05);
    \draw[fill = black] (5,0) circle(0.05);
    \draw[fill = black] (5.5,0) circle(0.05);
    \draw (7,0) [upup=1];
    \draw (9,0) -- (9,-1);
    \draw[fill = black] (9,0) circle(0.15);
    \draw (10,0) -- (10,-1);
    \draw[fill = black] (10,0) circle(0.15);
    \draw[fill = black] (11.5,0) circle(0.05);
    \draw[fill = black] (12,0) circle(0.05);
    \draw[fill = black] (12.5,0) circle(0.05);
    \draw[fill = black] (14,0) circle(0.15);
    \draw (14,0) -- (14,-1);
  \end{tikzpicture} \right) \not= 0$. In this case, we observe the following important equivalence which is intrinsic to $V_{n,p}$: 
  
  $$0 \not= a \left( \begin{tikzpicture}[scale = 0.2, baseline={(0,0.1)}]
    \draw (0,0) [upup=1];
    \draw (2,0) [upup=1];
    \draw[fill = black] (4.5,0) circle(0.05);
    \draw[fill = black] (5,0) circle(0.05);
    \draw[fill = black] (5.5,0) circle(0.05);
    \draw (7,0) [upup=1];
    \draw (9,0) -- (9,-1);
    \draw[fill = black] (9,0) circle(0.15);
    \draw (10,0) -- (10,-1);
    \draw[fill = black] (10,0) circle(0.15);
    \draw[fill = black] (11.5,0) circle(0.05);
    \draw[fill = black] (12,0) circle(0.05);
    \draw[fill = black] (12.5,0) circle(0.05);
    \draw[fill = black] (14,0) circle(0.15);
    \draw (14,0) -- (14,-1);
  \end{tikzpicture} \right) = b \left( \begin{tikzpicture}[scale = 0.2, baseline={(0,0.1)}]
    \draw (0,0) [upup=1];
    \draw (2,0) [upup=1];
    \draw[fill = black] (4.5,0) circle(0.05);
    \draw[fill = black] (5,0) circle(0.05);
    \draw[fill = black] (5.5,0) circle(0.05);
    \draw (7,0) [upup=1];
    \draw (9,0) -- (9,-1);
    \draw[fill = black] (9,0) circle(0.15);
    \draw (10,0) -- (10,-1);
    \draw[fill = black] (10,0) circle(0.15);
    \draw[fill = black] (11.5,0) circle(0.05);
    \draw[fill = black] (12,0) circle(0.05);
    \draw[fill = black] (12.5,0) circle(0.05);
    \draw[fill = black] (14,0) circle(0.15);
    \draw (14,0) -- (14,-1);
  \end{tikzpicture} \right) \iff a u_1 u_3 u_5 ... u_{2p-1} = b u_1 u_3 u_5 ... u_{2p-1}$$
  
 If we assume this equivalence to be true, note that we are done, because then we have that $a u_1 u_3 u_5 ... u_{2p-1} \alpha = b u_1 u_3 u_5 ... u_{2p-1} \alpha$, and since $\alpha \in \iota_* H_k(R_p)$, each of $u_1, u_3,u_5$ etc fix $\alpha$, and so we have that $a \alpha = b \alpha$. So, it remains to see why the above equivalence is true. 
 
 \
 
 Consider $a$. First, observe that $$a \left( \begin{tikzpicture}[scale = 0.2, baseline={(0,0.1)}]
    \draw (0,0) [upup=1];
    \draw (2,0) [upup=1];
    \draw[fill = black] (4.5,0) circle(0.05);
    \draw[fill = black] (5,0) circle(0.05);
    \draw[fill = black] (5.5,0) circle(0.05);
    \draw (7,0) [upup=1];
    \draw (9,0) -- (9,-1);
    \draw[fill = black] (9,0) circle(0.15);
    \draw (10,0) -- (10,-1);
    \draw[fill = black] (10,0) circle(0.15);
    \draw[fill = black] (11.5,0) circle(0.05);
    \draw[fill = black] (12,0) circle(0.05);
    \draw[fill = black] (12.5,0) circle(0.05);
    \draw[fill = black] (14,0) circle(0.15);
    \draw (14,0) -- (14,-1);
  \end{tikzpicture} \right) = a u_1 u_3 ... u_{2p-1} \left( \begin{tikzpicture}[scale = 0.2, baseline={(0,0.1)}]
    \draw (0,0) [upup=1];
    \draw (2,0) [upup=1];
    \draw[fill = black] (4.5,0) circle(0.05);
    \draw[fill = black] (5,0) circle(0.05);
    \draw[fill = black] (5.5,0) circle(0.05);
    \draw (7,0) [upup=1];
    \draw (9,0) -- (9,-1);
    \draw[fill = black] (9,0) circle(0.15);
    \draw (10,0) -- (10,-1);
    \draw[fill = black] (10,0) circle(0.15);
    \draw[fill = black] (11.5,0) circle(0.05);
    \draw[fill = black] (12,0) circle(0.05);
    \draw[fill = black] (12.5,0) circle(0.05);
    \draw[fill = black] (14,0) circle(0.15);
    \draw (14,0) -- (14,-1);
  \end{tikzpicture} \right)$$ 
  As a shorthand, write $a'  = a u_1 u_3 ... u_{2p-1}$. Write $a' = v w^*$, as in Definition \ref{vw*}. Since we have added $u_1 u_3 ... u_{2p-1}$ to the end of $a'$, $w$ must have cups at the first $p$ positions. Moreover, since we have assumed that $a \left( \begin{tikzpicture}[scale = 0.2, baseline={(0,0.1)}]
    \draw (0,0) [upup=1];
    \draw (2,0) [upup=1];
    \draw[fill = black] (4.5,0) circle(0.05);
    \draw[fill = black] (5,0) circle(0.05);
    \draw[fill = black] (5.5,0) circle(0.05);
    \draw (7,0) [upup=1];
    \draw (9,0) -- (9,-1);
    \draw[fill = black] (9,0) circle(0.15);
    \draw (10,0) -- (10,-1);
    \draw[fill = black] (10,0) circle(0.15);
    \draw[fill = black] (11.5,0) circle(0.05);
    \draw[fill = black] (12,0) circle(0.05);
    \draw[fill = black] (12.5,0) circle(0.05);
    \draw[fill = black] (14,0) circle(0.15);
    \draw (14,0) -- (14,-1);
  \end{tikzpicture} \right) \not= 0$, $w$ can have no further cups. Therefore, we must have that $w = \left( \begin{tikzpicture}[scale = 0.2, baseline={(0,0.1)}]
    \draw (0,0) [upup=1];
    \draw (2,0) [upup=1];
    \draw[fill = black] (4.5,0) circle(0.05);
    \draw[fill = black] (5,0) circle(0.05);
    \draw[fill = black] (5.5,0) circle(0.05);
    \draw (7,0) [upup=1];
    \draw (9,0) -- (9,-1);
    \draw[fill = black] (9,0) circle(0.15);
    \draw (10,0) -- (10,-1);
    \draw[fill = black] (10,0) circle(0.15);
    \draw[fill = black] (11.5,0) circle(0.05);
    \draw[fill = black] (12,0) circle(0.05);
    \draw[fill = black] (12.5,0) circle(0.05);
    \draw[fill = black] (14,0) circle(0.15);
    \draw (14,0) -- (14,-1);
  \end{tikzpicture} \right)$. Next, observe that since the last $n-2p$ loose strands are not affected by $w$, we must have that $v = a' \left( \begin{tikzpicture}[scale = 0.2, baseline={(0,0.1)}]
    \draw (0,0) [upup=1];
    \draw (2,0) [upup=1];
    \draw[fill = black] (4.5,0) circle(0.05);
    \draw[fill = black] (5,0) circle(0.05);
    \draw[fill = black] (5.5,0) circle(0.05);
    \draw (7,0) [upup=1];
    \draw (9,0) -- (9,-1);
    \draw[fill = black] (9,0) circle(0.15);
    \draw (10,0) -- (10,-1);
    \draw[fill = black] (10,0) circle(0.15);
    \draw[fill = black] (11.5,0) circle(0.05);
    \draw[fill = black] (12,0) circle(0.05);
    \draw[fill = black] (12.5,0) circle(0.05);
    \draw[fill = black] (14,0) circle(0.15);
    \draw (14,0) -- (14,-1);
  \end{tikzpicture} \right)$. Thus, we have that 
  
  \begin{equation*}
  \begin{split} 
  a' &= [a' \left( \begin{tikzpicture}[scale = 0.2, baseline={(0,0.1)}]
    \draw (0,0) [upup=1];
    \draw (2,0) [upup=1];
    \draw[fill = black] (4.5,0) circle(0.05);
    \draw[fill = black] (5,0) circle(0.05);
    \draw[fill = black] (5.5,0) circle(0.05);
    \draw (7,0) [upup=1];
    \draw (9,0) -- (9,-1);
    \draw[fill = black] (9,0) circle(0.15);
    \draw (10,0) -- (10,-1);
    \draw[fill = black] (10,0) circle(0.15);
    \draw[fill = black] (11.5,0) circle(0.05);
    \draw[fill = black] (12,0) circle(0.05);
    \draw[fill = black] (12.5,0) circle(0.05);
    \draw[fill = black] (14,0) circle(0.15);
    \draw (14,0) -- (14,-1);
  \end{tikzpicture} \right)] \cdot \left( \begin{tikzpicture}[scale = 0.2, baseline={(0,0.1)}]
    \draw (0,0) [upup=1];
    \draw (2,0) [upup=1];
    \draw[fill = black] (4.5,0) circle(0.05);
    \draw[fill = black] (5,0) circle(0.05);
    \draw[fill = black] (5.5,0) circle(0.05);
    \draw (7,0) [upup=1];
    \draw (9,0) -- (9,-1);
    \draw[fill = black] (9,0) circle(0.15);
    \draw (10,0) -- (10,-1);
    \draw[fill = black] (10,0) circle(0.15);
    \draw[fill = black] (11.5,0) circle(0.05);
    \draw[fill = black] (12,0) circle(0.05);
    \draw[fill = black] (12.5,0) circle(0.05);
    \draw[fill = black] (14,0) circle(0.15);
    \draw (14,0) -- (14,-1);
  \end{tikzpicture} \right)^* \\
  & = [b' \left( \begin{tikzpicture}[scale = 0.2, baseline={(0,0.1)}]
    \draw (0,0) [upup=1];
    \draw (2,0) [upup=1];
    \draw[fill = black] (4.5,0) circle(0.05);
    \draw[fill = black] (5,0) circle(0.05);
    \draw[fill = black] (5.5,0) circle(0.05);
    \draw (7,0) [upup=1];
    \draw (9,0) -- (9,-1);
    \draw[fill = black] (9,0) circle(0.15);
    \draw (10,0) -- (10,-1);
    \draw[fill = black] (10,0) circle(0.15);
    \draw[fill = black] (11.5,0) circle(0.05);
    \draw[fill = black] (12,0) circle(0.05);
    \draw[fill = black] (12.5,0) circle(0.05);
    \draw[fill = black] (14,0) circle(0.15);
    \draw (14,0) -- (14,-1);
  \end{tikzpicture} \right)] \cdot \left( \begin{tikzpicture}[scale = 0.2, baseline={(0,0.1)}]
    \draw (0,0) [upup=1];
    \draw (2,0) [upup=1];
    \draw[fill = black] (4.5,0) circle(0.05);
    \draw[fill = black] (5,0) circle(0.05);
    \draw[fill = black] (5.5,0) circle(0.05);
    \draw (7,0) [upup=1];
    \draw (9,0) -- (9,-1);
    \draw[fill = black] (9,0) circle(0.15);
    \draw (10,0) -- (10,-1);
    \draw[fill = black] (10,0) circle(0.15);
    \draw[fill = black] (11.5,0) circle(0.05);
    \draw[fill = black] (12,0) circle(0.05);
    \draw[fill = black] (12.5,0) circle(0.05);
    \draw[fill = black] (14,0) circle(0.15);
    \draw (14,0) -- (14,-1);
  \end{tikzpicture} \right)^* \\
  &= b' \\ 
  \end{split}
  \end{equation*}

\

For an example, $a' = b'$ might look like:

\

\begin{tikzpicture}[scale = 0.5]
\draw (1,4) [bupup=5];
\draw (2,4) [upup=3];
\draw (3,4) [supup=1];
\draw (8,4) [upup=1];

\draw (0,0) [dndn=1];
\draw (2,0) [dndn=1];
\draw (4,0) [dndn=1];
\draw (6,0) [dndn=1];

\draw (0,4) -- (0,2.5) -- (8,1.5) -- (8,0);
\draw (9,0) -- (9,1.5) -- (7,2.5) -- (7,4);
\draw[fill = black] (0,4) circle(0.029);
\draw[fill = black] (7,4) circle(0.029);
\draw[fill = black] (8,0) circle(0.029);
\draw[fill = black] (9,4) circle(0.029);
\end{tikzpicture}

\

Well definedness follows. 

\

\underline{Injectivity:}  We must show that if $v,w \in V_{n,p}$ are distinct link states then $h(v)$ and $h(w)$ are nonhomologous. If $v$ and $w$ do not have length one cups $\begin{tikzpicture}[scale = 0.2, baseline={(0,0.1)}]
    \draw (0,1) [upup=1];
  \end{tikzpicture}$ in the same locations, then the minimal intersections of $h(v)$ and $h(w)$ are different due to the observation we made in the ``Motivation:'' part of this proof, and thus $h(v) \not= h(w)$ by Corollary \ref{tlc}. Otherwise, $v$ and $w$ have their length one cups in the same locations. Then, since $v \not= w$, there is some length one cup such that the nesting around this cup in $v$ is more than in $w$. One then can multiply by an outermost (with respect to the nesting) $u_i$ to make $u_i v$ have length one cups at different positions as $u_i w$, and thus $h(u_i v) \not= h(u_i w)$, and thus $u_i h(v) \not= u_i h(w)$, and so $h(v) \not= h(w)$. The procedure for choosing $u_i$ is not mysterious, and is best illustrated by an example: 

\

$$v = \begin{tikzpicture}[scale = 0.2, baseline={(0,0.1)}]
    \draw (0,0) [bupup=5];
    \draw (1,0) [upup=3];
    \draw (2,0) [supup = 1];
    \draw (6,0) -- (6,-2);
    \draw[fill = black] (6,0) circle(0.15);
    \draw (7,0) [upup=1];
     \draw (9,0) -- (9,-2);
    \draw[fill = black] (9,0) circle(0.15);
    \end{tikzpicture}, w = \begin{tikzpicture}[scale = 0.2, baseline={(0,0.1)}]
    \draw (0,0) -- (0,-2);
    \draw[fill = black] (0,0) circle(0.15);
    \draw (5,0) -- (5,-2);
    \draw[fill = black] (5,0) circle(0.15);
    \draw (1,0) [upup=3];
    \draw (2,0) [supup = 1];
    \draw (6,0) [upup=3];
    \draw (7,0) [supup=1];
    \end{tikzpicture}$$
    
Observe that the length one cups of $v$ are both are positions $3$ and $8$. Thus, $h(v), h(w) \in \iota_* (A_3 \intersection A_8)$. But observe that the nesting around the length one cup at $3$ is larger for $v$ than $w$. Therefore, we choose $u_5$. Observe that 

$$u_5 v = \begin{tikzpicture}[scale = 0.2, baseline={(0,0.1)}]
    \draw (0,0) [upup=1];
    \draw (2,0) [upup = 1];
    \draw (4,0) [upup = 1];
    \draw (6,0) -- (6,-2);
    \draw[fill = black] (6,0) circle(0.15);
    \draw (7,0) [upup=1];
     \draw (9,0) -- (9,-2);
    \draw[fill = black] (9,0) circle(0.15);
    \end{tikzpicture}, u_5 w = \begin{tikzpicture}[scale = 0.2, baseline={(0,0.1)}]
    \draw (0,0) -- (0,-2);
    \draw (1,0) -- (1,-2);
    \draw[fill = black] (0,0) circle(0.15);
    \draw[fill = black] (1,0) circle(0.15);
    \draw (2,0) [upup=1];
    \draw (4,0) [upup = 1];
    \draw (6,0) [upup=3];
    \draw (7,0) [supup=1];
    \end{tikzpicture}$$
 
 \
 
 Notice therefore that $h (u_5 v) \in \iota_* H_k(A_1)$ but $h (u_5 w) \not\in \iota_* H_k(A_1)$. Thus, $u_5 h(v) = h(u_5 v) \not= h(u_5 w) = u_5 h(w)$ and thus, $h(v) \not= h(w)$.

\

\underline{Torsion:} Simply note that this module is just obtained by action on a set, and formal linearity. Therefore, torsion of $\alpha$ carries over verbatim to the coefficients of $V_{n,p}$. 
\end{proof}

\subsection{When invariants are trivial, the representation structure of $H_k(X)$ can be read off from its filtration $\scrF$} \label{ssrshxqrof}

\

Recall that in Definition \ref{fr}, we associated a filtration of retracts $\scrF$ to a given $\TL_n$-space. In this Subsection, we will present the proof of Theorem \ref{main}, which demonstrates how the representation structure of $H_k (X)$ can be read off from its filtration $\scrF$. 

\

Before we can state and prove Theorem \ref{main}, we need to introduce a few combinatorial notions which will appear in the statement of the theorem. Firstly, we need to introduce notation regarding compositions of an integer:

\

\begin{definition}[Notation for compositions] \label{nfc}
Given $m \in \N$, we will denote the collection of compositions of $m$ by $\Comp(m)$. If the reader needs reminding, a composition of $m$ is an ordered way to sum up to $m$. So, for example, $\Comp(3) = \{(1,1,1), (2,1), (1,2), (3)\}$. 

\

We will denote compositions by $\lambda$, since compositions can be viewed as generalizations of partitions, and hence pictorially speaking, compositions are generalizations of Young Diagrams. We will denote the number of rows of a given composition $\lambda$ by $\row(\lambda)$ - i.e $\row(\lambda)$ is the number of terms in the sum.
\end{definition}

\

One of the first results one encounters when learning the representation theory of Temperley-Lieb algebras is that $\dim V_{n,p} = \choose{n}{p} - \choose{n}{p-1}$ (See for example, \cite{tl1} (Ridout, Saint-Aubin), Section 2). This dimension is often denoted $d_{n,p}$, and we shall also denote it by $d_{n,p}$. Keeping this in mind, and motivated by the definitions of various families of symmetric functions, we arrive at the following definition: 

\

\begin{definition}[Defining $d^{r}_\lambda$] \label{ddrl}
Let $\lambda = (i_1,...,i_l)$ be a composition of $m$. For each $1 \le c \le l$, set $d^r_{\lambda,c} = d_{r- 2 \sum_{j=1}^{c-1} i_j, i_c}$. We then set

$$d^r_\lambda = \prod_{c = 1}^{l} d^r_{\lambda,c}$$
\end{definition}

\

We present an example for the reader to better understand the definition of $d^r_\lambda$. 

\

\begin{example}(Elucidating the definition of $d^r_\lambda$)
Consider $r = 8$. Noting that $\Comp(3) = \{(1,1,1), (2,1), (1,2), (3)\}$, we have:

\begin{equation*}
\begin{split} 
 d^8_{(1,1,1)}  &= d_{8,1} d_{6,1} d_{4,1} \\
 d^8_{(2,1)}  &= d_{8,2} d_{4,1} \\
 d^8_{(1,2)}  &= d_{8,1} d_{6,2} \\ 
 d^8_{(3)} &= d_{8,3} \\
 \end{split}
\end{equation*}
\end{example}

\

The above example would have made the definition of $d_\lambda$ clear to the reader. We now present a simple combinatorial Lemma which will help us prove our theorem. The Lemma describes how the various $d^r_\lambda$ are related.

\

\begin{lemma}[How different $d^r_\lambda$ are related] \label{hddr}
Let $\mu = (i_1,..,i_l) \in \Comp(m)$. Let $t \in \N$. Let $\lambda = (t, \mu)$ denote the composition of $m+t$ given by $(t,i_1,...,i_l)$. Then,

$$d^r_{(t, \mu)} = d_{r, t} \cdot d^{r - t}_\mu$$
\end{lemma}

\begin{proof}
We compute directly: 

\begin{equation*}
\begin{split}
d^r_{(t, \mu)} &= d_{r-0,t} \cdot \prod\limits_{c=1}^{\row(\mu)} d_{r - 2 (t + \sum_{j=1}^{c-1} i_j), i_c} \\ 
&= d_{r,t} \cdot \prod\limits_{c=1}^{\row(\mu)} d_{(r-2t) - 2 \sum_{j=1}^{c-1} i_j, i_c} \\ 
&= d_{r,t} \cdot d^{r-2t}_\mu\\ 
\end{split}
\end{equation*}

\end{proof}

\

We now prove the main theorem of this paper, Theorem \ref{main}, which gives an explicit formula for reading off the representation structure of homology from the associated filtration. 

\

\begin{theorem} \label{main}
Let $X$ be a surjective $\TL_n$-space and assume that $H_k(X)^{\textrm{TL}_n} = 0$. Let $\scrF = R_1 \supseteq R_2 \supseteq ... \supseteq  R_{\floor{\frac{n}{2}}}$ be the associated filtration. Let $\k$ be a field containing $\Q$. Let $\overset{\sim}{\TL_n} = \TL_n - \{1\}$ denote the ideal of $\TL_n$ containing all elements having cups. Then, for each $k$, and setting $r_p = n - 2p$ for each $p$, the $\TL_n$-representation structure of the subspace of homology which is the image of $\overset{\sim}{\TL_n} $ is given by:

\begin{equation*}
\overset{\sim}{\TL_n} H_k(X, \k) = \Ds\limits_{p=1}^{\floor{\frac{n}{2}}} V_{n,p}^{\Ds \left[\dim(H_k(R_p)) + \sum\limits_{p < q \le \floor{\frac{n}{2}}} \left( \sum\limits_{\lambda \in \Comp(q-p)} (-1)^{\row(\lambda)} \cdot d^{r_p} _\lambda \cdot \dim(H_k(R_q))\right) \right]}
\end{equation*}
\end{theorem}

\begin{proof}
We first note that $\overset{\sim}{\TL_n} H_k(X, \k)  = \textrm{Span}\{ \alpha \in H_k(X) \sep \alpha \in \iota_* H_k(A_j) \textrm{ for some j} \}$. Indeed,  given any $\beta \in H_k(X)$, note that $u_j \beta \in \iota_* H_k(A_j)$ for each $j$, and since $\overset{\sim}{\TL_n}$ is generated by $\{u_j\}_{j=1}^{n-1}$, we therefore have that $\overset{\sim}{\TL_n} H_k(X) \subseteq  \textrm{Span}\{ \alpha \in H_k(X) \sep \alpha \in \iota_* H_k(A_j) \textrm{ for some j} \}$. Conversely, if $\alpha \in \iota_* H_k(A_j)$, then $u_j \alpha = \alpha$, and therefore $\alpha \in \overset{\sim}{\TL_n} H_k(X, \k)$.

\

Next, observe that if $A_{i_1} \intersection ... \intersection A_{i_p}$ is any $p$ intersection, $\TL_n \cdot \iota_* H_k(A_{i_1} \intersection ... \intersection A_{i_p}) = \TL_n \cdot \iota_* H_k(R_p)$, since there are isomorphisms given by elements of $\TL_n$ which take each one to the other (by Lemma \ref{isch}). Therefore, we have that $\overset{\sim}{\TL_n} H_k(X)  = \Union_{p=1}^{\floor{\frac{n}{2}}} \TL_n \cdot \iota_* H_k(R_p)$. Of course, the terms in the above union have nontrivial redundancy. We must understand this redundancy, for which we may consider the following recursive procedure: 

\

\begin{itemize}
\item The base case is when $p = \floor{\frac{n}{2}}$, let $0 \not= \alpha \in \iota_* H_k(R_p)$.  In this case, since $\{1,3,5, ... ,2 \floor{\frac{n}{2}}-1 \}$ is a maximal subset of $\{1,...,n-1\}$ such that $|i_k - i_l| \ge 2$ for all elements $i_k, i_l$. Therefore, every nonzero $\alpha \in \iota_* H_k(R_{\floor{\frac{n}{2}}}, \k)$ satisfies the hypothesis of Theorem \ref{cmt}. Therefore, $R_{\floor{\frac{n}{2}}}$ contributes a term $V_{n,\floor{\frac{n}{2}}}^{\ds \dim H_k R_{\floor{\frac{n}{2}}}}$.
\item For the recursion, we suppose that, for any $q = p+1, p+2, .... , \floor{\frac{n}{2}}$, $R_q$ contributes a term $V_{n,q}^{\ds s_q}$. The set $\{1,2,3,...,s_q\}$ corresponds to a linearly independent set of homology classes $\alpha_1,...,\alpha_{s_q}  \in \iota_* H_k(R_q)$, but not all homology classes in $\iota_* H_k(R_q)$ live in the $\TL_n$-representation generated by these, because some classes belong to $V_{n,q'}$ for $q' > q$. We assume for the recursion that \\ $\iota_* H_k(R_q) = \Span\{\alpha_1,...,\alpha_{s_q}\} \ds (\Union_{q' > q} \TL_n \cdot \iota_* H_k(R_{q'}) ) \intersection \iota_* H_k(R_q)$ \\ 
\item We now do the inductive step. Extend a basis for $\left(\Union_{q > p} \TL_n \cdot \iota_* H_k(R_q)\right) \intersection \iota_* H_k(R_p) \subseteq \iota_* H_k(R_p)$ to a basis of $H_k(R_p)$, by appending a collection of linearly independent vectors $\{\alpha_1,...,\alpha_{s_p}\}$. Observe that in our inductive step, we assumed that for any $q > p$, $H_k(R_q) = \Span\{\alpha_1,...,\alpha_{s_q}\} \ds \left(\Union_{q' > q} \TL_n \cdot \iota_* H_k(R_{q'}) \right) \intersection \iota_* H_k(R_q)$, and thus all homology classes in $\iota_* H_k(R_{q})$ for $q > p$ have already been accounted for. Thus, must have that each $\alpha_i$ satisfies the hypothesis of Theorem \ref{cmt}, that is, if $T \supsetneq \{1,3,5, ... ,p\}$, $a \in \TL_n$, and $a \cdot \alpha_i \in \Intersection_{j \in T} \iota_* H_k(A_j)$, then in fact $a \cdot \alpha_i = 0$. Thus, by Theorem \ref{cmt}, we see that $\TL_n \cdot \alpha_i \isom V_{n,p}$. We deduce that $R_p$ contributes a term of $V_{n,p}^{\ds s_p}$ and moreover $\iota_* H_k(R_p) = \Span\{\alpha_1,...,\alpha_{s_p}\} \ds (\Union_{q > p} \TL_n \cdot \iota_* H_k(R_{q}) ) \intersection \iota_* H_k(R_p)$, thereby proving the inductive step. \\
\item Now, we notice that the ``redundancy'' comes from the fact that some of the homology classes in $\iota_* H_k(R_p)$ belong to copies of $V_{n,q}$ for $q > p$. Take $\alpha \in \iota_* H_k(R_q)$ for $q > p$. The elements in $\TL_n \cdot \alpha$ which reside in $\iota_* H_k(R_p)$ correspond to precisely those link states which have the first $p$ cups fixed in their positions, and the rest is allowed to vary. The ``rest'' amounts to placing $q-p$ cups in $n-2p$ locations, and so the number of link states with this property amounts choosing a link state of $V_{n-2p, q-p}$. Recall that $\dim V_{n-2p, q-p} = (\choose{n-2p}{q-p} - \choose{n-2p}{q-p-1})$, and we denote it by $d_{n-2p,q-p}$. We therefore conclude that $\dim(\TL_n \cdot \alpha \intersection \iota_* H_k(R_p)) = d_{n-2p,q-p}$. 
\item Therefore, $s_p$, the dimension of the span of all homology classes which are not obtained as part of some $V_{n,q}$ for $q > p$, is given by the recursion: 
\begin{equation*}
\begin{split}
&s_{\floor{\frac{n}{2}}} = \dim H_k (R_{\floor{\frac{n}{2}}}) \\ 
&s_p = \dim(H_k(R_p)) - \sum_{q = p+1}^{\floor{\frac{n}{2}} } d_{n-2p,q-p} \cdot s_q \\
\end{split}
\end{equation*}
\end{itemize}

\

What remains to be shown is that the above recursive formula yields the closed-form combinatorial formula which is in the statement of the theorem. For $p = \floor{\frac{n}{2}}$, observe that both formulae yield $\dim(H_k(R_{\floor{\frac{n}{2}}}))$, since, in both cases, the sum over $q > p$ is void.

\

Our inductive hypothesis is that for any $q > p$, $$s_q = \dim(H_k(R_q)) + \sum\limits_{q < q' \le \floor{\frac{n}{2}}} \left( \sum\limits_{\lambda \in \Comp(q'-q)} (-1)^{\row(\lambda)} \cdot d^{r_q} _\lambda \cdot \dim(H_k(R_{q'}))\right)$$

\

So, then:

\begin{equation*}
\begin{split}
s_p &= \dim(H_k(R_p)) - \sum_{q = p+1}^{\floor{\frac{n}{2}} } d_{n-2p,q-p} \cdot s_q \\ 
&= \dim(H_k(R_p)) - \sum_{q = p+1}^{\floor{\frac{n}{2}} } d_{n-2p,q-p} \cdot \left[ \dim(H_k(R_q)) + \sum\limits_{q < q' \le \floor{\frac{n}{2}}} \left( \sum\limits_{\mu \in \Comp(q'-q)} (-1)^{\row(\mu)} \cdot d^{r_q} _\mu \cdot \dim(H_k(R_{q'}))\right) \right] \\
\end{split}
\end{equation*}

\

\underline{Handling the first sum:}

\

Now, for $\lambda = (q-p) \in \Comp(q-p)$, observe that: \\
(1) $d^{r_p}_\lambda = d_{r_p,q-p} = d_{n-2p,q-p}$ \\
(2) $\row(\lambda) = 1$ and so $(-1)^{\row(\lambda)} = -1$ \\ 

\

Thus, we may suggestively rewrite 

\begin{equation*}
- \sum\limits_{q = p+1}^{\floor{\frac{n}{2}} } d_{n-2p,q-p} \cdot \dim(H_k(R_q)) = \sum\limits_{p < q \le \floor{\frac{n}{2}}} \left( \sum\limits_{\lambda \in \{(q-p)\} \subset \Comp(q-p)} (-1)^{\row(\lambda)} d^{r_p}_\lambda \dim(H_k(R_q)) \right)
\end{equation*}

\

and since relabeling $q$ by $q'$ in the above expression is harmless (since $q'$ does not appear), we have that: 

\begin{equation*}
-\sum\limits_{p < q' \le \floor{\frac{n}{2}}} d_{n-2p,q-p} \cdot \dim(H_k(R_q)) = \sum\limits_{p < q' \le \floor{\frac{n}{2}}} \left( \sum\limits_{\lambda \in \{(q'-p)\} \subset \Comp(q'-p)} (-1)^{\row(\lambda)} d^{r_p}_\lambda \dim(H_k(R_{q'})) \right)
\end{equation*}

\

\underline{Handling the second sum:}

\

Observe that given $q > p$ and $\mu \in \Comp(q'-q)$, $(q-p, \mu) \in \Comp(q'-p)$. Moreover, \\
(1) all such compositions are distinct, since the $\mu$s are distinct \\
(2) every composition in $\Comp(q'-p)$ is either of this form or is the single sum composition $(q'-p)$. \\

\

In other words, for each $q' > p$, $\Comp(q'-p) = \{(q'-p)\} \union \Union_{p < q < q'} \Comp(q' - q)$. 

\

Finally, observe that by Lemma \ref{hddr}, we have that $d^{r_p}_{(q-p, \mu)} = d_{r_p, q-p} \cdot d^{r_p - 2(q-p)}_\mu$, and observe that $r_p - 2(q-p) = n - 2p - 2(q-p) = n - 2q = r_q$, and thus we conclude that 

$$d^{r_p}_{(q-p, \mu)} =  d_{r_p, q-p} \cdot d^{r_q}_\mu$$

\

Therefore, combining what we have observed and noting that $\row(q-p, \mu) = 1+ \row(\mu)$, we see that: 

\

\begin{equation*}
{\tiny
\begin{split}
- \sum\limits_{q = p+1}^{\floor{\frac{n}{2}} } d_{n-2p,q-p} \sum\limits_{q < q' \le \frac{n}{2}} \left( \sum\limits_{\mu \in \Comp(q'-q)} (-1)^{\row(\mu)} \cdot d^{r_q} _\mu \cdot \dim(H_k(R_{q'}))\right) &=
- \sum\limits_{p < q < q' \le \floor{\frac{n}{2}}} \left( \sum\limits_{\mu \in \Comp(q'-q)} (-1)^{\row(\mu)} \cdot d_{n-2p,q-p} \cdot d^{r_q} _\mu \cdot \dim(H_k(R_{q'}))\right)  \\
 &= \sum\limits_{p < q < q' \le \floor{\frac{n}{2}}} \left( \sum\limits_{\mu \in \Comp(q'-q)} (-1)^{\row(\mu)+1} \cdot d^{r_p}_{(q-p, \mu)} \cdot \dim(H_k(R_{q'}))\right) \\
&= \sum\limits_{p < q' \le \floor{\frac{n}{2}}} \left( \sum\limits_{\lambda \in \Comp(q'-p) - \{(q'-p)\}} (-1)^{\row(\lambda)} \cdot d^{r_p}_\lambda \cdot \dim(H_k(R_{q'}))\right) \\
\end{split}
}
\end{equation*}

\

\underline{Putting together the handling of the first and second sums:} 

\

We saw that the first sum produced all required summands for compositions in the singletons $\{(q'-p)\}$ for all $q' > p$. We saw that the second sum produced all required summands for compositions except singletons, i.e $\lambda \in \Comp(q'-p) - \{(q'-p)\}$ for all $q' > p$. Putting these two together, we see that: 

\begin{equation*}
\begin{split}
s_p =& \dim(H_k(R_p)) + \textrm{(first sum)} + \textrm{(second sum)} \\
=& \dim(H_k(R_p)) + \sum\limits_{p < q' \le \floor{\frac{n}{2}}} \left( \sum\limits_{\lambda \in \{(q'-p)\}} (-1)^{\row(\lambda)} d^{r_p}_\lambda \dim(H_k(R_{q'})) \right) \\
&+ \sum\limits_{p < q' \le \floor{\frac{n}{2}}} \left( \sum\limits_{\lambda \in \Comp(q'-p) - \{(q'-p)\}} (-1)^{\row(\lambda)} \cdot d^{r_p}_\lambda \cdot \dim(H_k(R_{q'}))\right) \\ 
&= \dim(H_k(R_p)) + \sum\limits_{p < q \le \floor{\frac{n}{2}}} \left( \sum\limits_{\lambda \in \Comp(q-p)} (-1)^{\row(\lambda)} \cdot d^{r_p} _\lambda \cdot \dim(H_k(R_q))\right)
\end{split}
\end{equation*} 

\

We have therefore proved the inductive step, and this concludes the proof of the theorem.
\end{proof} 

\

\section{Topological stability} \label{sts}

\

The goal of this Section is to define a notion of topological stability of a chain of $\TL_n$-spaces, and observe that if $\{X_n\}_{n \ge N}$ is a topologically stable chain of $\TL_n$-spaces, then (an appropriate quotient of) their homology groups form a finitely generated $\LS$-module. This is philosophically interesting, because it is an analogue to a theorem in the foundational paper on representation stability (i.e \cite{rs1} (Church, Ellenberg, Farb), Section 6), which says that the homology groups of configuration spaces are finitely generated $\FI$-modules. 

\

\subsection{Defining topological stability} 

\

In the analogous story for symmetric groups, configuration spaces are natural candidates for a notion of topological stability. However, even there, it is not clear (at least to us), that there could not exist any other natural notion of topological stability of $S_n$-spaces, which has the property that chains of homology groups become finitely generated $\FI$-modules. Similarly, we do not claim here that the notion of topological stability that we present here is somehow overarching - indeed, there could be families of topological stable spaces which we have entirely missed. However, we will define a notion which appears natural to us: In Subsection \ref{ssfr}, we saw that associated to every topological action we may define an intrinsic filtration $\scrF$; recall that $\scrF$ is intrinsic in the sense that all $p$-intersections are homeomorphic by Lemma \ref{isch}. Since filtrations have a natural notion of stability, we may carry over this notion to define stability for a chain of $\TL_n$-spaces, as follows:

\

\begin{definition}[$p$-filtration stability] \label{fs}
Let $X_{N} \subseteq X_{N+1} \subseteq X_{N + 2} \subseteq .. $ be a chain of surjective $\TL_n$-spaces (each $X_n$ is a $\TL_n$-space), such that: \\
(1) The inclusions $X_n \subseteq X_{n+1}$ are inclusions of $\TL_n$-spaces (i.e $\TL_n$-equivariant). \\
(2) $H_k(X_n)^{\TL_n} = 0$ for all $k,n \ge 1$ (invariants are trival) \\

\

Let 
$$\scrF_n = R^{(n)}_1 \supseteq R^{(n)}_2 \supseteq ...  \supseteq R^{(n)}_{\floor{\frac{n}{2}}}$$

be the filtration associated to $X_n$. The chain $\{X_n\}_{n \ge N}$ is said to be \underline{$p$-filtration stable} if, for all $q \ge p$ and all $n \ge N$ $R^{(n)}_q = R^{(N)}_q$.
\end{definition}

\

\subsection{$p$-filtration stability results in representation stability of homology groups}

\

We will show in this Subsection how $p$-filtration stability implies the representation stability of certain quotients of images of homology groups. These quotients will be quotients in the Grothendieck group, and thus we require some notation:

\

\begin{definition}[Quotients in the Grothendieck group] \label{qgg}
Suppose that $U_1,..,U_l$ are indecomposable representations of an algebra $A$ with the property that representations in the split Grothendieck group generated by  $U_1,..,U_l$ have unique direct sum decompositions in terms of $U_1,..,U_l$. Suppose that $V = \Ds_{i=1}^l U_i^{\ds m_i}$. Then, we denote:

$$\quotient{V}{[U_{i_1}], ... , [U_{i_d}]} = \Ds_{i \in \{1,...,l\} - \{i_1,...,i_d\} } U_i^{\ds m_i}$$

We call it a quotient in the Grothendieck group generated by $U_1,..,U_l$, since the Grothendieck bracket of the resulting representation is equal to $\quotient{[V]}{[U_{i_1}], ... , [U_{i_d}]}$ in the quotient of the split Grothendieck group by $[U_{i_1}], ... , [U_{i_d}]$.
\end{definition}

\

We can now state and prove the theorem of this Subsection, which is the second result stated in the introduction, and which can be viewed as a corollary of the work we did in Section \ref{shgtlr}:

\

\begin{corollary} \label{fsirs}
Let $\{X_n\}_{n \ge N}$ be a $p$-filtration stable chain of $\TL$-spaces. Let $\k$ be a field containing $\Q$. Let $\overset{\sim}{\TL_n} = \TL_n - \{1\}$ denote the ideal of $\TL_n$ containing all elements having cups. Then, for each $k$, $\{\quotient{\overset{\sim}{\TL_n} H_k(X_n, \k)}{[V_{n,1}], ... , [V_{n,p-1}]}\}_{n \ge N}$ is a finitely generated $\LS$-module.
\end{corollary}

\begin{proof}
Since $\{X_n\}_{n \ge N}$ is $p$-filtration stable, we know that for any $n \ge N$ and $q \ge p$, $R^{(n)}_q = R^{(N)}_q$. Theorem \ref{main} gives an explicit decomposition of $\overset{\sim}{\TL_n} H_k(X, \k)$ in terms of the filtration $\scrF$. Furthermore, observe that since $R_q = R_p$ for $q \ge p$ for $n \ge N$, and since the contribution of $V_{n,p}$ in our formula only depends on $R_q$ for $q \ge p$, we see that $$\{\quotient{\overset{\sim}{\TL_n}  H_k(X_n, \k)}{[V_{n,1}], ... , [V_{n,p-1}]}\}_{n \ge N} = \{V_{n,p}^{\ds s_1} \ds .... \ds V_{n, p_{max}}^{\ds s_{p_{max}}} \}_{n \ge N}$$ where $p_{max}$ is the largest $q$ such that $R_q \not= \{*\}$, which is the same for all $n \ge N$ by filtration stability. Moreover, note that for $m,n \ge N$ the induced inclusion on homology takes, in the notation of the proof of Theorem \ref{cmt}, $h\left(\begin{tikzpicture}[scale = 0.2, baseline={(0,0.1)}]
    \draw (0,0) [upup=1];
    \draw (2,0) [upup=1];
    \draw[fill = black] (4.5,0) circle(0.05);
    \draw[fill = black] (5,0) circle(0.05);
    \draw[fill = black] (5.5,0) circle(0.05);
    \draw (7,0) [upup=1];
    \draw (9,0) -- (9,-1);
    \draw[fill = black] (9,0) circle(0.15);
    \draw[fill = black] (10,0) circle(0.05);
    \draw[fill = black] (10.5,0) circle(0.05);
    \draw[fill = black] (11,0) circle(0.05);
    \draw[fill = black] (12,0) circle(0.15);
    \draw (12,0) -- (12,-1);
  \end{tikzpicture}\right)$ to $h\left(\begin{tikzpicture}[scale = 0.2, baseline={(0,0.1)}]
    \draw (0,0) [upup=1];
    \draw (2,0) [upup=1];
    \draw[fill = black] (4.5,0) circle(0.05);
    \draw[fill = black] (5,0) circle(0.05);
    \draw[fill = black] (5.5,0) circle(0.05);
    \draw (7,0) [upup=1];
    \draw (9,0) -- (9,-1);
    \draw[fill = black] (9,0) circle(0.15);
    \draw (10,0) -- (10,-1);
    \draw[fill = black] (10,0) circle(0.15);
    \draw[fill = black] (11.5,0) circle(0.05);
    \draw[fill = black] (12,0) circle(0.05);
    \draw[fill = black] (12.5,0) circle(0.05);
    \draw[fill = black] (14,0) circle(0.15);
    \draw (14,0) -- (14,-1);
  \end{tikzpicture} \right)$, where the first link state has $m-2p$ loose strands while the second has $n-2p$ loose strands. The rest of the map is determined by cyclicity, and thus if we restrict the map to any standard representation we see that the restriction $(i_{m,n})_* \restriction_{V_{m,q}} V_{m,q} \to V_{n,q}$ is just the usual $\LS$-module map associated with the chain of standard representations $\{V_{n,q}\}_{n \in \N}$, in the sense of Corollary \ref{csrs}. Thus, in light of Observation \ref{sfn}, we conclude that $\{\quotient{\overset{\sim}{\TL_n} H_k(X_n, \k)}{[V_{n,1}], ... , [V_{n,p-1}]}\}_{n \ge N}$ is a finitely generated $\LS$-module.
\end{proof}

\

\subsection{Examples of topological stability} \label{ssets}

\

The goal of this Subsection is to give a couple of examples of topological stability. This Section will likely not only be useful to a reader who is interested in topological stability, but will also be useful to a reader who wants to understand topological actions as in Section \ref{sdtatla}.

\

\begin{example}\label{ex1}
Just as before, let $d_{n,p} = \dim V_{n,p} = \choose{n}{p} - \choose{n}{p-1}$. Let $n \ge 2p$. Let $X_{n,p}$ be the wedge of $d_{n,p}$ copies of $S^2$. Let the filtration associated to $\TL_n$-action on $X_{n,p}$ be

$$\scrF_{n,p} = \twedge_{d_{n-2, p-1}} S^2 \supseteq \twedge_{d_{n-4, p-2}} S^2 \supseteq ... \supseteq \twedge_{d_{n-2p, 1}} S^2 \supseteq S^2 \supset \{*\} \supset ... \supset \{*\}$$

In particular note that if $A_1,...,A_{n-1}$ are the retracts given by $u_1,...,u_{n-1}$, we have that $A_i \isom \twedge_{d_{n-2, p-1}} S^2 $ for each $i$. There is a well defined action of  $\TL_n$ on $X_{n,p}$ with this filtration, which we have not yet defined, and will do so after the example.  The reader is encouraged to define this action for some small values of $n$ and $p$, to appreciate the fact that a ``naive approach'' of writing down the actions from a topological perspective requires a little thought, since one needs to be careful that the retractions must be chosen to be compatible with one another. Nevertheless, one can do it. In this spirit, we have illustrated below an example for $n = 5$ and $p = 2$.

\

\begin{figure}[H]
\centering
\begin{tikzpicture}
   \begin{polaraxis}[grid=none, axis lines=none]
     \addplot[mark=none, color = black, very thick, domain=-18:18,samples=300] {-cos(5*x)};
   \end{polaraxis}
   \begin{polaraxis}[grid=none, axis lines=none]
     \addplot[mark=none, color = black, very thick, domain=18:54,samples=300] {-cos(5*x)};
   \end{polaraxis}
   \begin{polaraxis}[grid=none, axis lines=none]
     \addplot[mark=none, color = black, very thick, domain=54:90,samples=300] {-cos(5*x)};
   \end{polaraxis}
    \begin{polaraxis}[grid=none, axis lines=none]
     \addplot[mark=none, color = black, very thick, domain=90:126,samples=300] {-cos(5*x)};
   \end{polaraxis}
    \begin{polaraxis}[grid=none, axis lines=none]
     \addplot[mark=none, color = black, very thick, domain=126:162,samples=300] {-cos(5*x)};
   \end{polaraxis}
   \node at (0,3.4) {$\mathbf{3,1}$};
   \node at (2.8,6.8) {$\mathbf{1,4}$};
   \node at (2.8,0) {$\mathbf{3}$};
   \node at (6.2,5.6) {$\mathbf{4,2}$};
   \node at (6.3,1.4) {$\mathbf{2}$};
   
   \draw [->, thick, red] (5,4.5) -- (3,5); 
   \draw [->, thick, red] (5,2.5) -- (1.5,3.5);
   \draw [->, thick, red] (3,1.5) -- (1,3.5); 
   
   \draw [->, thick, blue] (2.8,5.5) -- (5.5,5); 
   \draw [->, thick, blue] (1.5,3.2) -- (5.5,2); 
   \draw [->, thick, blue] (2.8,1.3) -- (4.8,4.4); 
\end{tikzpicture}
\caption{Above is a schematic of $X_{5,2}$. This has a filtration $S^2 \twedge S^2 \supseteq S^2$. In particular, note that each $A_i \isom S^2 \twedge S^2$. The numbers next to the copies of $S^2$ indicate the minimal intersection containing that copy. For example, $1,4$ means that the minimal intersection containing that copy is $A_1 \intersection A_4$. The action is as follows: each $u_i$ preserves $A_i$, takes $A_{i-1}$ and $A_{i+1}$ isomorphically to $A_i$ and retracts $A_j$ to $A_i$ for $|j-i| \ge 2$. However, these isomorphisms and retractions must be chosen carefully so as to not lead to any contradictions - this is really the heart of the matter. In red arrows above, we have drawn the action of $u_1$. In blue arrows above, we have drawn the action of $u_2$. The action of $u_3$ is similar to the action of $u_2$ and the action of $u_4$ is similar to the action of $u_1$.}
\end{figure}

\

In case the reader is wondering how this action is defined and how each $u_i$ is defined in general, we note that we may ``work backwards'' to construct the action. That is, we can define this action by assigning one copy of $S^2$ for each $(n,p)$ link state, and $\TL_n$ action permutes the copies of $S^2$ via the actions on link states. Whenever $u_i$ of a link state is $0$, $u_i$ acts on the corresponding copy of $S^2$ by sending it to the wedge point $\{*\}  = \intersection_{i=1}^n A_i$.

\

In any case, for $p$ fixed, $\{X_{n,p}\}_{n \in \N}$ is $p$-filtration stable since the filtration stabilizes after position $\ge p$ and so for all $n \ge 1$, $H_n(X_{n,p})^{\TL_n} = 0$. Furthermore, since every copy of $S^2$ is fixed by some $u_i$, we have that $H_n(X_{n,p}) = \overset{\sim}{\TL_n} H_n(X_{n,p})$, where $\overset{\sim}{\TL_n} = \TL_n - \{1\}$ denotes the ideal of $\TL_n$ containing all elements having cups. We therefore see that $\quotient{H_k(X_{n,p})}{[V_{n,1}], ..., [V_{n,p-1}]}$ is a finitely generated $\LS$-module. Explicitly, we know that for $i \not= 2$, $H_i(X_{n,p}) = 0$ for all $n$ and $p$, and for $i=2$, as an $\LS$-module:

$$\{H_2(X_{n,p})\}_{n \ge 2p} = \{V_{n,p}\}_{n \ge 2p} =: V_{*,p}$$

At this point, we would like to emphasize the importance of Theorem \ref{main} in deducing the above result. This justifies us to use the ``backwards cheating trick'' to construct actions. Without the theorem, a priori one might have believed that there are other actions possible which would give us other representations. 
\end{example}

\

\begin{example}\label{ex2}
Let $\{c_n\}_{n \ge 2}$ be a sequence of integers such that $c_n \ge \max\{d_{n,2}, n-1\}$ for each $n$. Given the $c_n$-torus $T^{c_n} = S^1 \times S^1 \times ... \times S^1$, we may construct a non-surjective action of $\TL_n$ on $T^{c_n}$ by utilizing the projection maps onto the first $n-1$ coordinates (which we will explicitly define). We will then use this to construct a surjective action on a bigger space, and the sequence of such bigger spaces will be $2$-filtration stable, with filtration: 

$$T^{c_n} \supseteq S^1 \supseteq \{*\} \supseteq ... \supseteq \{*\}$$

\

Let us first, for concreteness, define the non-surjective action explicitly. We think of $T^{c_n}$ as $(\quotient{\R}{\Z})^{c_n}$. If we let $p_j$ denote the projection onto the $j$th coordinate, we define $r_i$ for $i \in \{1,2,3, ... , n-1\}$ by:  

$$p_j r_i(x_1,...,x_{c_n}) = \begin{cases} x_{i-1} + x_i + x_{i+1} \delta_{i+1 \le n-1} & \textrm{ if $ j = i$ } \\ 0  & \textrm{ if $ j \not= i$ }  \end{cases}$$

\

This is a well defined action since, if $|i_1 - i_2| \ge 2$ then $r_{i_1} r_{i_2} = 0 = r_{i_2} r_{i_1}$, and 

$$p_j r_i r_{i\pm1} r_i(x_1,...,x_{c_n}) = \begin{cases} x_{i-1} + x_i + x_{i+1} \delta_{i+1 \le n-1} & \textrm{ if $ j = i$ } \\ 0  & \textrm{ if $ j \not= i$ } \end{cases} = r_i(x_1,...,x_{c_n})$$

\

Out of this non-surjective action we construct a surjective action on another space as follows by taking $n-1$ copies of $T^{c_n}$ and glue them along the retracts. More precisely, let $(T^{c_n})_1, ... , (T^{c_n})_{n-1}$ be the $n-1$ copies of the tori, let $(T^{c_n})_{i j}$ denote the $j$th copy of $S^1$ in $(T^{c_n})_i$, that is: $0 \times 0 \times ... \times 0 \times (S^1)_j \times 0 \times ... \times 0_{c_n} \subset (T^{c_n})_i$. And we define

$$X_n = \quotient{ \Disj\limits_{i=1}^{n-1} (T^{c_n})_i}{(T^{c_n})_{kl} \sim (T^{c_n})_{lk} \textrm{ for $|k -l| \ge 2$, $k,l \le n-1$}}$$

We construct action on $X_n$ by letting $u_i$ be the retraction onto $(T^{c_n})_i$. More explicitly:

\

The $1$ skeleton of $X$ carries an action of $\TL_n$ and has a filtration given by (since one copy of far intersections are glued):

$$\scrF = \twedge_{n-1} S^1 \twedge_{c_n-(n-1)} S^1 \supseteq S^1 \supseteq \{*\} \supseteq ... \supseteq \{*\}$$

It is not hard to see that the one skeleton is a the wedge of two spaces, one with filtration $\scrF = \twedge_{n-1} S^1 \supseteq S^1 \supseteq \{*\} \supseteq ... \supseteq \{*\}$ and the other which is a wedge of copies of $V_{n,1}$. By an argument similar to Example \ref{ex1} by replacing $S^2$ with $S^1$, we see that the one skeleton decomposes as $V_{n,2} \ds V_{n,1}^{\ds c_n - d_{n,2}}$.

\

\underline{Formal construction of space:}  Let $X_{n,2}, X_{n,1}$ denote the analogous space of Example \ref{ex1} by replacing $S^2$ with $S^1$. Notice that the $1$ skeleton of the above space is equal to  $\left[X_{n,2} \twedge (\twedge_{(c_n - d_{n,2})} X_{n,1})\right]_i$. For each $i$, take the lengths of all loops in $\left[X_{n,2} \twedge (\twedge_{(c_n - d_{n,2})} X_{n,1})\right]_i$ to be $1$, and likewise take the lengths of all loops in the $1$-skeleton of $T^{c_n}_i$ to be $1$. Then, for each $i$, fix an isometry:

$$\phi_i: \left[X_{n,2} \twedge (\twedge_{(c_n - d_{n,2})} X_{n,1})\right]_i  \isom \textrm{ $1$-skeleton of $T^{c_n}_i$ which is given by considering components }$$

with the property that for each $i$, $j$, $$\phi_i(u_j \cdot \left[X_{n,2} \twedge (\twedge_{(c_n - d_{n,2})} X_{n,1})\right]_i) = T^{c_n}_{ij}$$

\

We then define the action of $\TL_n$ on $X_n$ by setting: For each $i$, 

\begin{equation*}
\begin{split}
u_i(x_1,...,x_{c_n}) &= \begin{cases} \sum\limits_{s=1}^{c_n} \phi_i \circ u_i \circ \phi_{i \pm 1}^{-1} x_s \in T^{c_n}_{i} & \textrm{ if $(x_1,...,x_{c_n}) \in (T^{c_n})_{i \pm 1}$} \\  \sum\limits_{s=1}^{c_n} \phi_j \circ u_i \circ \phi_{j}^{-1} x_s \in T^{c_n}_{ji}& \textrm{ if $(x_1,...,x_{c_n}) \in (T^{c_n})_{j}$ where $|j-i| \ge 2$ } \end{cases} \\ 
&= \begin{cases} \sum\limits_{s=1}^{c_n} \phi_i \circ u_i \circ \phi_{i \pm 1}^{-1} x_s \in T^{c_n}_{i} & \textrm{ if $(x_1,...,x_{c_n}) \in (T^{c_n})_{i \pm 1}$} \\ r_i(x_1,...,x_{c_n}) \in T^{c_n}_{ji} & \textrm{ if $(x_1,...,x_{c_n}) \in (T^{c_n})_{ji}$ where $|j-i| \ge 2$ } \end{cases} \\ 
\end{split}
\end{equation*}

\

Then, since the action of $\TL_n$ on $X_{n,2} \twedge (\twedge_{(c_n - d_{n,2})} X_{n,1})$ is well defined, the action of $\TL_n$ on $X_n$ is well defined, as is best seen by considering the first line in the above equation together with the reminder that $T^{c_n}_{ji} \sim T^{c_n}_{ij}$. Below is a schematic of $X_4$.

\

\begin{figure}[H]
\centering
\begin{tikzpicture}
\draw (0,0) -- (0.7,0.7);
\draw (0,0) -- (0,1);
\draw (0,0) -- (-0.7,0.7);
\draw (0,1) -- (0.7,1.7);
\draw (0,1) -- (-0.7,1.7);
\draw (0.7,0.7) -- (0.7,1.7);
\draw (-0.7,0.7) -- (-0.7,1.7);
\draw (0.7,1.7) -- (0,2.2);
\draw (-0.7,1.7) -- (0,2.2);

\draw (1.4,0) -- (2.1,0.7);
\draw (1.4,0) -- (1.4,1);
\draw (1.4,0) -- (0.7,0.7);
\draw (1.4,1) -- (2.1,1.7);
\draw (1.4,1) -- (0.7,1.7);
\draw (2.1,0.7) -- (2.1,1.7);
\draw (0.7,0.7) -- (0.7,1.7);
\draw (2.1,1.7) -- (1.4,2.2);
\draw (0.7,1.7) -- (1.4,2.2);

\draw (0.7,0.7) -- (0.7,-0.2);
\draw (0.7,0.7) -- (1.2,-0.3);
\draw (0.7,0.7) -- (0.2,-0.3);
\draw (0.2,-0.3) -- (0.2,-1);
\draw (1.2,-0.3) -- (1.2,-1);
\draw (0.7,-0.2) -- (0.2,-1);
\draw (0.7,-0.2) -- (1.2,-1);
\draw (0.7,-1.7) -- (0.2,-1);
\draw (0.7,-1.7) -- (1.2,-1);

\draw [->, thick, red] (1.5,1) -- (0.8,1); 
\draw [->, thick, red] (1.5,0.5) -- (0.8,0.8); 
\draw [->, thick, red] (1.5,1.5) -- (0.8,1.2); 
\draw [->, thick, red] (0.8,-0.8) -- (0.1,0.9); 

\draw [->, thick, blue] (-0.4,1) -- (0.6,-1.1); 
\draw [->, thick, blue] (1.7,1) -- (0.8,-1.1); 

\node at (-1,2) {$\mathbf{1}$};
\node at (2.2,2) {$\mathbf{3}$};
\node at (0,-1.5) {$\mathbf{2}$};
\end{tikzpicture}
\caption{A schematic of $X_4$. We have taken $c_4 = 3$. The cubes of course depict $3$-tori. The numbers $1,2,3$ depict the retracts $A_1, A_2, A_3$ respectively - each of which is a $3$-torus. The red arrows depict the action of $u_1$ - $A_3$ retracts to the copy of $S^1 \subset A_1 \intersection A_3$ while $A_2$ is taken isomorphically to $A_1$. The blue arrows depict the action of $u_2$ - both $A_1$ and $A_3$ are taken isomorphically to $A_2$. The action of $u_3$ is similar to the action of $u_1$.}
\end{figure} 

\

Now, since the filtrations $\{\scrF_n\}_{n \ge 2}$ stabilize at $p = 2$,  and since $H_k(X_n)^{\TL_n} = 0$ for each $n$, we know from Theorem \ref{fsirs} that for each $k$, $\quotient{\overset{\sim}{\TL_n} H_k(X_n)}{[V_{n,1}]}$ is a finitely generated $\LS$-module, where $\overset{\sim}{\TL_n} = \TL_n - \{1\}$ denotes the ideal of $\TL_n$ containing all elements having cups. Since $\overset{\sim}{\TL_n} H_k(X_n) = H_k(X_n)$ in this case, we have that $\quotient{H_k(X_n)}{[V_{n,1}]}$ is a finitely generated $\LS$-module. Moreover, by Theorem \ref{main}, we know that the generators of this $\LS$ module all lie in $R_2 = S^1$. It therefore follows that for $k >1$, $\quotient{H_k(X_n)}{[V_{n,1}]} = 0$ and for $k = 1$, $\quotient{H_k(X_n)}{[V_{n,1}]} = V_{n,2}$. Therefore, as an $\LS$-module,

$$\{\quotient{H_1(X_n)}{[V_{n,1}]}\}_{n \ge 2} = \{V_{n,2}\}_{n \ge 2} =: V_{*,2}$$

\

Notice that while $\{\quotient{H_1(X_n)}{[V_{n,1}]}\}_{n \ge 2} = \{V_{n,2}\}_{n \ge 2}$ is a finitely generated $\LS$-module, $\{H_1(X_n)\}_{n \ge 2}$ is NOT a finitely generated $\LS$-module in general since, by Theorem \ref{main}, we may explicitly decompose $H_k(X_n) = V_{n,2} \ds V_{n,1}^{\ds (c_n-d_{n,2})}$, and we may choose $c_n$ to grow faster than quadratically, so that $H_k(X_n)$ fails to be finitely generated. This demonstrates why we needed to take a quotient in the statement of Corollary \ref{fsirs}. Lastly, notice that we constructed this family of $\TL_n$-spaces by, to an extent, ``working backwards''. Theorem \ref{main} suggests to us that we ought to be able to work backwards to construct an arbitrary $\TL_n$-space satisfying $H_k(X)^{\TL_n} = 0$ for all $k$.
\end{example} 

\

\section{Future directions} \label{sfd}

\

(1) We view the results of this paper as an occurrence of a more general phenomenon, which we now outline: 
\begin{itemize}
\item Suppose that $\{A_n\}_{n \in \N}$ is a natural chain of algebras under inclusions, possibly depending on some parameters that lie in a field $\k$. \\
\item Specialize these parameters so that no relations involve ``$+$'' signs, and so that every element of $\k$ that appears in a relation is either $0$ or $1$. \\ 
\item $0$ and $1$ have topological meaning (null-homotopic and the identity, respectively), and therefore we can study topological actions of $A_n$. \\ 
\item Try to figure out how to decompose homology groups as $A_n$ representations (these are of course honest algebra representations, and one can make full use of addition to study them). \\ 
\item Develop a theory for representation stability for a chain $\{A_n\}_{n \in \N}$, either actions on finite sets as we did in this paper or by extracting sufficiently nice families of representations of $A(\infty)$ (see \cite{si1}, \cite{si2}, \cite{si4} (Okounkov, Lieberman, Vershik, Nessonov) for analogies). \\ 
\item Define a natural notion of topological stability based on the structure theory of $A_n$-spaces. \\ 
\item Deduce whether or not topological stability implies representation stability of homology groups. \\ 
\end{itemize}

We view this paper as a starting point, and aim to study a host of other important families of algebras along the lines of the procedure outlined above.

\

(2) One might try to apply the results of this paper to find group (co)homological obstructions to $\TL_n$ actions on groups (as a monoid), since a $\TL_n$ action on a group induces an action on its $K(G,1)$ space. Unfortunately though, this action need not be surjective, and therefore one would need to find a clever way to construct a surjective space out of the resulting space, and perhaps write the homology of the resulting space in terms of our original space. One would also need to account for the fact that $H_k(K(G,1))^{\TL_n}$ could be nontrivial. One might also want to compute induced action on the fundamental group of the new space one would construct. Example \ref{ex2} is an example of how to construct a new space with a surjective action from a space with a non-surjective action. Notice that, there, the original action we begin with is an action on a $c_n$-torus, which is the $K(G,1)$-space of $\Z^{c_n}$, and therefore the action that we begin with in Example \ref{ex2} corresponds to an action on the group $\Z^{c_n}$.

\

(3) A perhaps surprising fact that we learnt from Theorem \ref{main} is that, for any $\TL_n$-space $X$ with $H_k(X)^{\TL_n} = 0$, the only representations that occur in the nontrivial part as summands are the standard representations $V_{n,p}$. Our proof of this was constructive in nature. In \cite{tl2} (Belletete, Ridout, Saint-Aubin), the authors provide a classification of indecomposable modules over $\TL_n$ (It is not so straightforward, since we are are working at the root of unity $q = e^{\frac{i \pi}{3}})$. It is conceivable that another proof of the fact that we unearthed in this paper could be obtained by studying the structure theory of the other indecomposable on the list, and comparing it with structural obstructions of the condition $H_k(X)^{\TL_n} = 0$. While constructive proofs like we present in this paper are usually preferred, this other proof would also be valuable, since it would shed insight both into structural features present in topologically obtained representations and structural features of the various indecomposables.

\

\section{Acknowledgements} 

\

We would like to thank Mikhail Khovanov for his insightful comments and his kind encouragement. We would like to thank Sun Woo Park for helping us proofread the paper. We would like to thank the anonymous referees for their valuable comments, which helped improve the paper substantially and correct some of the oversights in early versions of this paper.

\

\bibliography{TATLRS}{}
\bibliographystyle{plain}

\end{document}